\documentclass{amsart}
\usepackage{amsmath,amsfonts}
\usepackage{amscd,amsthm}

\theoremstyle{plain} 

\newtheorem{theorem}{Theorem}[section]
\newtheorem{proposition}[theorem]{Proposition}
\newtheorem{corollary}[theorem]{Corollary}
\newtheorem{lemma}[theorem]{Lemma}
\newtheorem{remark}[theorem]{Remark}

\numberwithin{equation}{section}

\renewcommand{\theremark}{\theremark.\arabic{remark}}

\newcommand{\val}[2]{\mathop{ {\it Val}\left( #1 \right)\left( #2 \right) }}

\newcommand{\vall}[3]{\mathop{ {\it Val}_{ #3 }\!\left( #1 \right)\left( #2 \right) }}
\newcommand{\forest}[1]{\mathop{ \mathcal{T}_{#1}}}
\newcommand{\R}{\mathop{\mathbb{R}}}
\newcommand{\Q}{\mathop{\mathbb{Q}}}
\newcommand{\Z}{\mathop{\mathbb{Z}}}
\newcommand{\N}{\mathop{\mathbb{N}}}
\newcommand{\C}{\mathop{\mathbb{C}}}
\newcommand{\bigo}[1]{\mathop{\mathcal{O}(#1)}}
\newcommand{\class}[1]{\mathop{\mathcal{C}_{(#1)}}}

\newcommand{\diag}[1]{\mathop{\it diag(#1_1,\dots,#1_n)}}

\title[Lagrange's inversion formula. Non--Analytic Linearization Problems.]{The Lagrange inversion formula on non--Archimedean fields. Non--Analytical Form of Differential and Finite Difference Equations.}

\author[Timoteo Carletti]{}

\email{carletti@math.unifi.it}



\begin{document}
\maketitle

\setcounter{page}{1}

\centerline{\scshape  Timoteo Carletti}
\medskip

{\footnotesize
\centerline{ Dipartimento di Matematica "U. Dini"}
\centerline{ viale Morgagni 67/A }
\centerline{ 50134 Firenze, Italy }
}
\medskip


\bigskip
\begin{quote}{\normalfont\fontsize{8}{10}\selectfont
{\bfseries Abstract.}
The classical Lagrange inversion formula is extended to analytic and non--analytic inversion problems
on non--Archimedean  fields. We give some applications to the field
of formal Laurent series in $n$ variables, where the non--analytic inversion formula gives  explicit formal solutions of general 
semilinear differential and $q$--difference equations. 

We will be interested in linearization problems for germs of diffeomorphisms (Siegel center problem) and vector fields. In addition to analytic results, we give sufficient condition for the linearization to belong to some Classes of ultradifferentiable germs, closed under composition and derivation, including Gevrey Classes. We prove that Bruno's condition is sufficient for the linearization to belong to the same Class of the germ, whereas new conditions weaker than Bruno's one are introduced if one allows the linearization to be less regular than the germ. This generalizes to dimension $n> 1$ some results of~\cite{CarlettiMarmi}. Our formulation of the Lagrange inversion formula by mean of trees, allows us to point out the strong similarities existing between the two linearization problems, formulated (essentially) with the same functional equation. For analytic vector fields of $\C^2$ we prove a quantitative estimate of a previous qualitative result of~\cite{MatteiMoussu} and we compare it with a result of~\cite{YoccozPerezMarco}.
\par}
\end{quote}

\section{Introduction}

Let $k$ be a field of characteristic zero complete
with respect to a non--trivial absolute value $|\; |$ and let $k'$
denote its residue field. When $k=\mathbb{R}$ or $\mathbb{C}$, the
classical Lagrange inversion formula
(see~\cite{Lagrange}\footnote{J. H. Lambert was the first interested in determining the roots $x$ of the equation $x^m+px=q$ developing it in an infinite series~\cite{Cajori}(but also see~\cite{Lubet}). His results stimulated J. Lagrange, who first generalized the method to solve the equation $a-x+\phi (x)=0$ for an analytic function $\phi$, and then he applied the idea to the Kepler problem: solving the elliptic motion of a point mass planet
about a fixed point, according to the law of the inverse square.
To do this Lagrange studied the possibility of inverting the fundamental relation between the mean
anomaly, $M$, and the eccentric anomaly, $E$: $E=e\sin E + M$, being $e$ the eccentricity of the orbit.},
\cite{Dieudonne} chapter VIII, section 7 or~\cite{Sansone} p. 286, for the $1$--dimensional case, and to~\cite{Good} for the multidimensional case) says that if $G$ is an analytic function in a neighborhood
of $w\in k$ then there exists a unique solution $h=H(u,w)$ of
\begin{equation}
h =u G \left( h \right) + w \, ,
 \label{inversion}
\end{equation}
provided that $|u|$ is sufficiently small. The solution $h=H(u,w)$ depends
analytically on $u$ and $w$ and its Taylor series with respect to $u$
is explicitly given by the formula:
\begin{equation}
H\left(u,w\right)=w+\sum_{n \geq 1}\frac{u^n}{n!}
\frac{d^{n-1}}{dw^{n-1}}(G\left( w \right))^n \, .
\label{eq:lagrangeorig}
\end{equation}
\indent In sections~\ref{sec:anallag}, after recalling some elementary notions of theory of
analytic functions on non--Archimedean fields, we give two
generalizations of~\eqref{inversion}:
 in the $n$--dimensional vector space $k^n$, $k$ {\em non--Archimedean} when $G$ is an analytic function (Corollary~\ref{existence}), and
 for non--analytic $G$  (Theorem~\ref{nonanex}). To deal with this second
case we rewrite the Lagrange inversion formula by means of the tree formalism. We refer to~\cite{Gessel} and references therein for a combinatorial proof of the Lagrange inversion formula using the tree formalism.

\indent In sections~\ref{sec:someappl} and~\ref{sec:nonanvf} 
we will give some applications of the
previous results in the setting of the formal Laurent series with applications to some dynamical systems problems. The idea of using trees in non--linear small divisors problems (in particular Hamiltonian) is due to H. Eliasson~\cite{Eliasson} who introduced trees in his study of the absolute convergence of Lindstedt series. The idea has been further developed by many authors (see, for example, ~\cite{ChierchiaFalcolini,Gallavotti1,Gallavotti2} always in the context of Hamiltonian KAM theory, see also~\cite{BerrettiGentile1} which we take as reference for many definitions concerning trees). The fact that these formulas should be obtained by a suitable generalization of Lagrange's inversion formula was first remarked by Vittot~\cite{Vittot}.

When $k$ is the field of formal Laurent series ${\mathbb C}((z))$, we
consider the vector space: $\mathbb{C}^n\left(\left(z_1,\ldots,z_n\right)\right)$; the
non-analytic inversion problem can be applied to obtain the
solution of semilinear differential or $q$--difference equations
in an explicit (i.e. not recursive) form. Our results are
formulated so as to include general first--order $U$--differential
semilinear equations~\cite{Du} and semilinear convolution equations.
In particular we will study (section~\ref{sec:someappl}) the Siegel center problem~\cite{Herman,Yoccoz} for
 analytic and non--analytic germs of $(\mathbb{C}^n,0)$, $n \geq 1$, and (section~\ref{sec:nonanvf}) the Problem of linearization of analytic~\cite{Bruno} 
and non--analytic vector fields of $\C^n$, $n \geq 1$. 
The reader interested only in the Siegel center problem may find useful 
to assume Proposition~\ref{thecoefficients}  and to skip the reading of 
the whole of sections~\ref{sec:anallag} and~\ref{sec:proofs}. The same is 
true for those interested in the linearization of vector fields, assuming Proposition~\ref{thecoefficientsvf} and reading the rest of section~\ref{sec:nonanvf}, 
even if they will find several useful definitions in section~\ref{sec:someappl}.

In~\cite{CarlettiMarmi} authors began the study of the Siegel center problem in some ultradifferentiable algebras of  $\mathbb{C}\left(\left(z\right)\right)$, 
here we generalize these results to dimension $n\geq 1$.

Consider two Classes of formal power series $\mathcal{C}_1$ and $\mathcal{C}_2$ of $\C^n\left[\left[z_1,\dots,z_n\right]\right]$ closed with respect to the composition and derivation. For example the Class of germs of analytic functions of $(\C^n,0)$ or Gevrey--$s$ Classes, $s>0$ (i.e. series ${\bf F}=\sum_{\alpha\in\N^n}{\bf f}_{\alpha}z^{\alpha}$ for which there exist $c_1,c_2>0$ such that $|{\bf f}_{\alpha}|\leq c_1 c_2^{|\alpha|}(|\alpha|!)^s$, for all $\alpha\in\N^n$). Let $A\in GL(n,\C)$ and ${\bf F}\in \mathcal{C}_1$ such that ${\bf F}({\bf z})=A{\bf z}+\dots$, we say that ${\bf F}$ is {\em linearizable} in $\mathcal{C}_2$ if there exists ${\bf H}\in \mathcal{C}_2$, tangent to the identity, such that:
\begin{equation*}
{\bf F}\circ {\bf H}({\bf z})={\bf H}(A{\bf z})
\end{equation*}
When $A$ is in the Poincar\'e domain (see \S~\ref{subsec:soemknownres}), the results of Poincar\'e~\cite{Poincare1} and Koenigs~\cite{Koenigs} assure that ${\bf F}$ is linearizable in $\mathcal{C}_2$. When $A$ is in the Siegel domain (see \S~\ref{subsec:soemknownres}), the problem is harder, the only trivial case is $\mathcal{C}_2=\C^n\left[\left[z_1,\dots,z_n\right]\right]$ ({\em formal linearization}) for which one only needs to assume $A$ to be {\em non--resonant}.

In the analytic case we recover the results of Bruno~\cite{Bruno} and R\"ussmann~\cite{Russmann}, whereas in the non--analytic case new arithmetical conditions are introduced (Theorem~\ref{maintheorem}). Consider the general case where both $\mathcal{C}_1$ and $\mathcal{C}_2$ are different from the Class of germs of analytic function of $(\C^n,0)$, if one requires $\mathcal{C}_1=\mathcal{C}_2$, once again the Bruno condition is sufficient, otherwise if $\mathcal{C}_1\subset\mathcal{C}_2$ one finds new arithmetical conditions, weaker than the Bruno one.

In section~\ref{sec:nonanvf} we will consider the following differential equation:
\begin{equation}
\dot {\bf z}=\frac{d{\bf z}}{dt}={\bf F}({\bf z}) \, ,
\label{leq:diffeqa}
\end{equation}
where $t$ is the time variable and ${\bf F}$ is a formal power series in the $n\geq 1$  variables $z_1,\dots , z_n$, with coefficients in $\C^n$, without constant term: ${\bf F}=\sum_{\alpha \in \N^n,|\alpha|\geq 1}{\bf F}_{\alpha}z^{\alpha}$, and we are interested in the behavior of the solutions near the singular point ${\bf z}=0$. 

A basic but clever idea has been introduced by Poincar\'e (1879), which consists in reducing the system~\eqref{leq:diffeqa} with an appropriate change of variables, to a simpler form: the {\em normal form}. In~\cite{Bruno} several results are presented in the analytic case (namely ${\bf F}$ is a convergent power series). Here we generalize such kind of results to the case of non--analytic ${\bf F}$, with a {\em diagonal, non--resonant linear part}. More precisely considering the same Classes of formal power series as we did for the Siegel Center Problem , we take an element ${\bf F}\in \mathcal{C}_1$ with a diagonal, non--resonant linear part, $A{\bf z}$, and we look for sufficient conditions on $A$ to ensure the existence of a change of variables ${\bf H}\in \mathcal{C}_2$ (the {\em linearization}), such that in the new variables the vector field reduces to its linear part. We will show that the Bruno condition is sufficient to linearize in the same class of the given vector field, whereas in the general case, $\mathcal{C}_1\subset \mathcal{C}_2$, new arithmetical conditions, weaker than the Bruno one, are introduced (Theorem~\ref{the:mainvf}). Finally in the case of analytic vector field of $\mathbb{C}^2$, the use of the continued fraction and of a best description of the accumulation of small divisors (due to the Davie counting function~\cite{Davie}), allows us to improve (Theorem~\ref{the:mainvfn2}) the results of Theorem~\ref{the:mainvf}, giving rise to (we conjecture) an optimal estimate concerning the domain of analyticity of the linearization. This gives a quantitative estimate of some previous results of~\cite{MatteiMoussu} and~\cite{YoccozPerezMarco}.

In our formulation we emphasize the strong similarities existing between this problem and the Siegel Center Problem, which becomes essentially the same problem; in fact once we reduced each problem to a Lagrange inversion formula (on some appropriate setting) we get the same functional equation to solve.

\section{The Lagrange inversion formula on non--Archimedean fields}
\label{sec:anallag}

In this section we generalize the Lagrange inversion  formula for
analytic and non--analytic functions on complete, ultrametric fields
of characteristic zero. In the first part we give for completeness
some basic definitions and properties of non--Archimedean fields,
referring to Appendix~\ref{ultrametricappendix} and
to~\cite{Serre,Christol,Chenciner} for a more detailed discussion.
We end the section introducing some elementary facts concerning trees.

\subsection{Statement of the Problem}
\label{subsect:ultramatric}

\indent Let $(k,|\;|)$ be a non--Archimedean field~\footnote{The 
reader can keep in mind the following two main models of 
non--Archimedean fields: the formal Laurent series and the $p$--adic
 numbers (examples a) and b) page 65 of~\cite{Serre}).} of
characteristic zero, where $|\;|$ is a ultrametric absolute value
: $|x+y|\le \sup (|x|,|y|)$ for all $x,y\in k$. Moreover we assume
that $k$ is complete and the norm is {\em non--trivial}. Let $a$
be a real number such that $0<a<1$, given any $x \in k$ we define
the real number $ v \left( x \right)$ by: $\lvert x \rvert = a^{v
\left( x \right)}$, the {\em valuation}\footnote{From the
properties of $|\;|$ it follows that the valuation satisfies, for
all $x,y \in k$ : $v \left( x \right) = + \infty$ if and only if
$x=0$; $v \left( xy \right) = v \left( x \right)+v \left( y
\right)$; $v \left( 1 \right) = 0$; $v \left( x+y \right) \geq
\inf \left( v \left( x \right),v \left( y \right)\right)$.} of
$x$.

\indent Since $k$ is non--Archimedean one has the following
elementary but  fundamental result:

\begin{proposition}
\label{converges}
Let $\left( x_n \right)_{n \in \mathbb{N}}$ be a sequence with $x_n
\in k$. Then $\sum x_n$ converges if and only if $x_n \rightarrow 0$.
\end{proposition}

\indent Let $n \in \mathbb{N}$, we introduce the
$n$--dimensional vector space $k^n$ and using the ultrametric
absolute value, defined on $k$, we introduce a {\em norm} $|| \cdot ||:k^n \rightarrow \mathbb{R}_+$
\begin{equation}
\label{defnorm}
|| x || = \sup_{1\leq i \leq n}| x_i| \quad x=(x_1,\ldots, x_n) \in k^n \, ,
\end{equation}
which results an ultrametric one and verifies a Schwartz--like
inequality, for $x,y \in k^n$ then: $|x\cdot y| \leq ||x||\, ||y||$, where $x\cdot y=\sum_{i=1}^n x_iy_i$ is a scalar product.

\indent This norm induces a topology, where the open balls are
defined\footnote{One could define~\cite{Serre} the open polydisks: $P_0(x,\rho)=\left\{y\in k^n:\forall i,1\leq i\leq n:|x_i-y_i|<\rho_i\right\}$, for some $x\in k^n$ and $\rho\in\mathbb{R}_+^n$. Clearly the induced
topology is equivalent to the previously defined one.} by
\begin{equation}
\label{defopenset}
B_0(x,r)=\left\{ y\in k^n:||x-y||<r\right\}
\end{equation}
for $x\in k^n$ and $r\in \mathbb{R}_+$. We will denote the closed ball
with $B(x,r)=\overline{B}_0(x,r)$.

\indent Let $r>0$ and let us consider a function $G: B_0(0,r)\subset k^n \rightarrow k^{n \times l}$, i.e. for $x \in B_0(0,r)$ and for all $1 \leq i \leq n ,1 \leq j \leq l$:
\begin{equation*}
G \left( x \right) = \left( G_{ij}\left( x \right)  \right)_{ij},
\quad G_{ij}\left( x \right) \in k.
\end{equation*}
Given $w\in k^n$, $u \in k^l$ and $G$ as above, we consider the
following problem:
\begin{quotation}
Solve with respect to $h\in k^n$, the {\em multidimensional
non--analytic Lagrange inversion problem}:
\begin{equation}
h = \Lambda \left[w + G\left( h \right) \cdot u\right] ,\label{multilagrange}
\end{equation}
where $\Lambda$ is a $k^n$--additive, $k^\prime$--linear, 
non--expanding operator (i.e. $\lvert \lvert \Lambda w \rvert
\rvert \leq \lvert \lvert w \rvert \rvert$ for all $w\in k^n$).
\end{quotation}

We will prove the existence of a solution of~\eqref{multilagrange}
using \emph{trees}. We will now recall some elementary facts concerning trees; we refer to~\cite{Harary} for a more
complete description.

\subsection{The Tree formalism} \label{sssec:treeform}

A tree is a connected acyclic graph, composed by {\em nodes} and
{\em lines}  connecting together two or more nodes. Among trees we
consider {\em rooted trees}, namely trees with an extra node, not
included in the set of nodes of the tree, called the {\em earth},
and an extra line connecting the earth to the tree, the {\em root
line}. We will call {\em root} the only node to which the earth is
linked. The existence of the root introduces a {\em partial
ordering} in the tree: given any two nodes\footnote{To denote
nodes we will use letters: $u,v,w,\ldots$, with possible
sub-indices. Lines will be denoted by $\ell$, the line exiting
from the node $u$ will be denoted by $\ell_u$.} $v,v^{\prime}$,
then $v\geq v^{\prime }$ if the (only) path connecting the root
$v_1$ with $v^{\prime}$, contains $v$. The {\em order} of a tree is the
number of its nodes. The {\em forest} $\mathcal{T}_N$ is the disjoint
union of all trees~\footnote{Here we consider only
\emph{semitopological} trees (see~\cite{BerrettiGentile2}), we
refer to~\cite{Gallavotti1} for the definition
of \emph{topological trees}.} with the same order $N$.

The {\em degree of a node}, $\deg v$, is the number of  incident lines
with the node. Let $m_v = \deg v-1$, that is the number of lines
entering into the node $v$ w.r.t. the partial ordering, if $m_v=0$ we
 will say that $v$ is an {\em end node}; for the
root $v_{1}$, because the root line doesn't belong to the lines of
the tree, we define $m_{v_1}=\deg v_1$, in this way $m_{v_1}$ also
represents the number of lines entering in the root. Let $\vartheta$ be a rooted
tree, for any $v \in \vartheta$ we denote by $L_v$ the 
{\em set of lines entering into $v$}; if $v$ is an end 
node we will set $L_v=\emptyset$.

Given a rooted tree $\vartheta$ of order $N$, we  can view it as
the union of its root and the subtrees $\vartheta^i$ obtained from
$\vartheta$ by detaching the root. Let $v_1$ be the root of
$\vartheta$ and $t=m_{v_1}$, we define the {\em standard
decomposition} of $\vartheta$ as: $\vartheta = \{ v_1 \} \cup
\vartheta^1 \cup \dots \cup \vartheta^t$, where $\vartheta^i \in
\mathcal{T}_{N_i}$ with $N_1+ \ldots +N_t=N-1$.

Using the definition of $m_v$ we can associate  uniquely to a
rooted tree of order $N$ a vector of $\mathbb{N}^N$, whose
components are just $m_v$ with $v$ in the tree~\cite{Vittot}. Thus
$\mathcal{T}_N =\{ (m_1,\ldots,m_N) \in \mathbb{N}^N :
\sum_{i=1}^N m_i=N-1, \quad \sum_{i=j}^N m_i \leq N-j \quad
\forall j=1, \cdots ,N \}$. We can then rewrite the standard
decomposition of $\vartheta$ as: $\vartheta = \left(
t,\vartheta^1, \ldots , \vartheta^t \right)$ where the subtrees
satisfy: $\vartheta^i \in \mathcal{T}_{N_i}$ with $N_1+ \ldots
+N_t=N-1$.

\begin{center}
  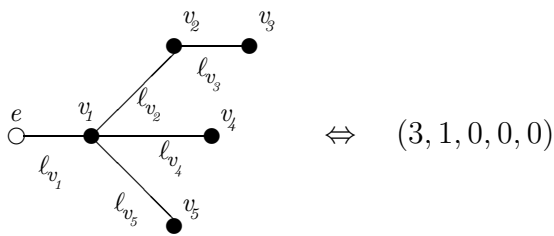
\begin{figure}[ht]
\setlength{\unitlength}{4144sp}%
\begingroup\makeatletter\ifx\SetFigFont\undefined%
\gdef\SetFigFont#1#2#3#4#5{%
  \reset@font\fontsize{#1}{#2pt}%
  \fontfamily{#3}\fontseries{#4}\fontshape{#5}%
  \selectfont}%
\fi\endgroup%
\begin{picture}(1898,1343)(128,-519)
\put(2026, 29){\makebox(0,0)[lb]{\smash{\SetFigFont{12}{14.4}{\familydefault}{\mddefault}{\updefault}$\Leftrightarrow \quad (3,1,0,0,0)$}}}
\put(406,-241){\makebox(0,0)[lb]{\smash{\SetFigFont{5}{6.0}{\rmdefault}{\mddefault}{\itdefault}1}}}
\put(361,-196){\makebox(0,0)[lb]{\smash{\SetFigFont{8}{9.6}{\rmdefault}{\mddefault}{\itdefault}v}}}
\put(901,254){\makebox(0,0)[lb]{\smash{\SetFigFont{10}{12.0}{\rmdefault}{\mddefault}{\itdefault}$\ell$}}}
\put(946,209){\makebox(0,0)[lb]{\smash{\SetFigFont{8}{9.6}{\rmdefault}{\mddefault}{\itdefault}v}}}
\put(991,164){\makebox(0,0)[lb]{\smash{\SetFigFont{5}{6.0}{\rmdefault}{\mddefault}{\itdefault}2}}}
\put(1171,-376){\makebox(0,0)[lb]{\smash{\SetFigFont{10}{12.0}{\rmdefault}{\mddefault}{\itdefault}v}}}
\put(1216,-421){\makebox(0,0)[lb]{\smash{\SetFigFont{6}{7.2}{\rmdefault}{\mddefault}{\itdefault}5}}}
\put(1261,434){\makebox(0,0)[lb]{\smash{\SetFigFont{10}{12.0}{\rmdefault}{\mddefault}{\itdefault}$\ell$}}}
\put(1306,389){\makebox(0,0)[lb]{\smash{\SetFigFont{8}{9.6}{\rmdefault}{\mddefault}{\itdefault}v}}}
\put(1351,344){\makebox(0,0)[lb]{\smash{\SetFigFont{5}{6.0}{\rmdefault}{\mddefault}{\itdefault}3}}}
\put(1171,749){\makebox(0,0)[lb]{\smash{\SetFigFont{10}{12.0}{\rmdefault}{\mddefault}{\itdefault}v}}}
\put(1216,704){\makebox(0,0)[lb]{\smash{\SetFigFont{6}{7.2}{\rmdefault}{\mddefault}{\itdefault}2}}}
\put(1621,749){\makebox(0,0)[lb]{\smash{\SetFigFont{10}{12.0}{\rmdefault}{\mddefault}{\itdefault}v}}}
\put(1666,704){\makebox(0,0)[lb]{\smash{\SetFigFont{6}{7.2}{\rmdefault}{\mddefault}{\itdefault}3}}}
\put(541,209){\makebox(0,0)[lb]{\smash{\SetFigFont{10}{12.0}{\rmdefault}{\mddefault}{\itdefault}v}}}
\put(586,164){\makebox(0,0)[lb]{\smash{\SetFigFont{6}{7.2}{\rmdefault}{\mddefault}{\itdefault}1}}}
\put(766,-376){\makebox(0,0)[lb]{\smash{\SetFigFont{10}{12.0}{\rmdefault}{\mddefault}{\itdefault}$\ell$}}}
\put(811,-421){\makebox(0,0)[lb]{\smash{\SetFigFont{8}{9.6}{\rmdefault}{\mddefault}{\itdefault}v}}}
\put(856,-466){\makebox(0,0)[lb]{\smash{\SetFigFont{5}{6.0}{\rmdefault}{\mddefault}{\itdefault}5}}}
\put(1396,164){\makebox(0,0)[lb]{\smash{\SetFigFont{10}{12.0}{\rmdefault}{\mddefault}{\itdefault}v}}}
\put(1441,119){\makebox(0,0)[lb]{\smash{\SetFigFont{6}{7.2}{\rmdefault}{\mddefault}{\itdefault}4}}}
\put(1036,-61){\makebox(0,0)[lb]{\smash{\SetFigFont{10}{12.0}{\rmdefault}{\mddefault}{\itdefault}$\ell$}}}
\put(1081,-106){\makebox(0,0)[lb]{\smash{\SetFigFont{8}{9.6}{\rmdefault}{\mddefault}{\itdefault}v}}}
\put(1126,-151){\makebox(0,0)[lb]{\smash{\SetFigFont{5}{6.0}{\rmdefault}{\mddefault}{\itdefault}4}}}
\thinlines
\put(1126,614){\circle*{90}}
\put(1576,614){\circle*{90}}
\put(181, 74){\circle{90}}
\put(631, 74){\circle*{90}}
\put(1126,-466){\circle*{90}}
\put(1351, 74){\circle*{90}}
\put(631, 74){\line( 1, 1){517.500}}
\put(1126,614){\line( 1, 0){450}}
\put(631, 74){\line(-1, 0){405}}
\put(629, 71){\line( 1,-1){517}}
\put(631, 74){\line( 1, 0){720}}
\put(136,164){\makebox(0,0)[lb]{\smash{\SetFigFont{10}{12.0}{\rmdefault}{\mddefault}{\itdefault}e}}}
\put(316,-151){\makebox(0,0)[lb]{\smash{\SetFigFont{10}{12.0}{\rmdefault}{\mddefault}{\itdefault}$\ell$}}}
\end{picture}
    \caption{A rooted tree of order $5$, $v_i$, $i=1,\ldots,5$, are the nodes whereas $\ell_{v_i}$, $i=1,\ldots,5$, are the lines. The earth is denoted by the letter $e$, $\ell_{v_1}$ is the root line and $v_3,v_4,v_5$ are end nodes. On the right we show the standard decomposition of this tree.}
    \label{figtree5}
  \end{figure}
\end{center}

\indent In the following we will also use {\em labeled rooted
trees}. A labeled rooted tree of order $N$ is an element of
$\mathcal{T}_N$ together with $N$ labels: $\alpha_1,\ldots, \alpha_N$. We
can think that the label $\alpha_i$ is attached to the i--th node of
the standard decomposition of the tree. The label is nothing else
that a function from the set of nodes of a tree to some set,
usually a subset of $\mathbb{Z}^m$ for some integer $m$. When
needed we denote a labeled rooted tree of order $N$ with the
couple $(\vartheta,\mathbf{\alpha})$, where $\vartheta \in
\mathcal{T}_N$ and $\mathbf{\alpha}=(\alpha_1,\ldots,\alpha_N)$ is the vector
label.

\subsection{The non--analytic Lagrange inversion formula}

We are now able to extend equation~\eqref{inversion}, the
classical analytic Lagrange inversion formula, to the setting of
paragraph~\ref{subsect:ultramatric}. We refer the reader to Appendix~\ref{ultrametricappendix} for a brief introduction to the theory of analytic functions on $k^n$ (norms, Cauchy estimates, etc ... ).

Let $N\in \mathbb{N}^*$, $U$ and $V$ be open subsets of,
respectively, $k^l$ and $k^n$, and $\vartheta \in \mathcal{T}_N$. We
define the function ${\it V\!al}_{\Lambda}:\forest{N} \times U \times V \ni (\vartheta,u,w)\mapsto \vall{\vartheta}{u,w}{\Lambda}\in k^n$ as follows
\begin{equation}
  \label{defvallambda}
\vall{\vartheta}{u,w}{\Lambda}=\begin{cases} \Lambda \left(
G\left( \Lambda w \right) \cdot u \right)& \text{if
$\vartheta\in\mathcal{T}_1$}\\
      \frac{1}{t!}\Lambda \left[ d^{t}G\left( \Lambda w \right) \left(
\vall{\vartheta^1}{u, w}{\Lambda}, \ldots , \vall{\vartheta^t}{u,
w}{\Lambda}  \right) \cdot u \right] & \text{otherwise}
\end{cases}
\end{equation}
where $\vartheta=(t,\vartheta^1,\ldots,\vartheta^t)$ is the
standard decomposition of the tree and $\Lambda :k^n\rightarrow k^n$.

\begin{remark}
\label{rem:different} For $t\geq 1$ and $v^{\left( 1 \right)},
\dots,v^{\left( t \right)}\in k^n$ we recall that
\begin{equation*}
d^{t}G_{ij}\left( w \right) \left( v^{\left( 1 \right)}, \dots,
v^{\left( t \right)} \right)=\sum_{l_1,\dots,l_t=1}^n D_{l_1}\dots
D_{l_t} G_{ij}\left( w \right) v^{\left( 1 \right)}_{l_1} \dots
v^{\left( t \right)}_{l_t},
\end{equation*}
with $v^{(i)}=(v_1^{(i)},\ldots,v_n^{(i)})$ and $D_{l_i}$ are the
$l_i$--th partial derivatives of $G$ at $w$ (see
Appendix~\ref{ultrametricappendix}).
\end{remark}

\indent We can then state the following existence Theorem:

\begin{theorem}[Non--analytic case]
  \label{nonanex}
Let $n,l$ be positive integers. Let $u\in k^l$, $w \in k^n$ and
for $1 \leq i \leq n$ let $G_i=\left( G_{i1}, \dots
,G_{il}\right)$, with $G_{i}\in k^n \left[ \left[ X_1, \dots, X_n
\right] \right]$. Let $\Lambda:k^n\rightarrow k^n$ be a
$k^n$--additive, $k^{\prime}$--linear and non--expanding operator.
Assume that for $\left( r_1, \dots ,r_n \right)\in \R_+^n$ the
$G_i$'s are convergent in $B\left( 0,r_i \right)$,
 set\footnote{See Appendix\ref{ultrametricappendix} for the definition  $||\cdot||_{r}$.} $M_i=||G_i||_{r_i}>0$, $r=\min_i r_i$ and $M=\max_i M_i$.
Then equation~\eqref{multilagrange} has the unique solution
\begin{equation}
\label{nonansol} H\left( u,w \right)=\Lambda w+ \sum_{N\geq
1}\sum_{\vartheta \in \mathcal{T}_{N}}
\vall{\vartheta}{u,w}{\Lambda}\, ,
\end{equation}
where ${\it Val}_{\Lambda}: \mathcal{T}_N\times B_{0}(0,r/M)\times
B_{0}(0,r) \rightarrow B_{0}(0,r)$, has been defined in
(\ref{defvallambda}). Moreover for any fixed $\vartheta$ the function ${\it
Val}_{\Lambda}$ is continuous, series~\eqref{nonansol}
converges on $B_{0}(0,r/M)\times B_{0}(0,r)$ and the map
$(u,w)\mapsto H(u,w)\in B_{0}(0,r)$ is continuous.
\end{theorem}

\begin{remark}
Since $\Lambda$ is not $k^n$--linear, ${\it Val}_{\Lambda}$ and
$H$ cannot be analytic. However the non--expanding condition
implies that $\Lambda$ is Lipschitz continuous from which the
regularity properties of ${\it Val}_{\Lambda}$ and $H$ follow.

\indent If $\Lambda$ is $k^n$--linear then we have the following Corollary (particular
case of the previous Theorem) with $w^\prime$ instead of $\Lambda w$ and $G^\prime\cdot u^\prime$ 
instead of $\Lambda (G\cdot u)$, which extends the analytic Lagrange inversion formula~\eqref{inversion}.
\end{remark}

\begin{corollary}[Analytic case]
\label{existence} Let $n,l$ be positive integers. Let $u^\prime\in k^l$,
$w^\prime \in k^n$ and for $1 \leq i \leq n$ let $G^\prime_i=\left( G^\prime_{i1},
\dots ,G^\prime_{il}\right)$, with $G^\prime_{i}\in k^n \left[ \left[ X_1,
\dots, X_n \right] \right]$. Assume that for $\left( r_1, \dots
,r_n \right)\in \R_+^n$ the $G^\prime_i$'s
 are convergent in $B\left( 0,r_i \right)$, let $M_i=||G^\prime_i||_{r_i}>0$, $r=\min_i r_i$ and $M=\max_i M_i$. Then
\begin{equation}
\label{solution} H\left( u^\prime,w^\prime \right)=w^\prime+ \sum_{N\geq
1}\sum_{\vartheta \in \mathcal{T}_{N}}  \val{\vartheta}{u^\prime,w^\prime}\, ,
\end{equation}
is the unique solution of~\eqref{multilagrange} with $\Lambda (G\cdot u)=G^{\prime}\cdot u^\prime$ and $\Lambda w=w^{\prime}$. The function $\val{\vartheta}{u^\prime,w^\prime}$
 is nothing else that the function $\vall{\vartheta}{u,w}{\Lambda}$ with
$\Lambda (G \cdot u)= G^\prime \cdot u^\prime$ and $\Lambda w=w^\prime$. Moreover $H$  is
analytic in $B_0\left( 0,\frac{r}{M} \right) \times B_0\left( 0,r
\right)$.
\end{corollary}

\section{Proofs.}
\label{sec:proofs}

This section is devoted to the proof of Theorem~\ref{nonanex} and Corollary~\ref{existence}.

\noindent{\em Proof of Theorem~\ref{nonanex}.} 
Using the fact that $\Lambda$ is non--expanding the
uniqueness of the solution can be proved easily. Let $H_1$ and
$H_2$ be two solutions of~\eqref{multilagrange}, then
\begin{equation*}
  ||H_1-H_2||=||\Lambda \left[ \left( G\left( H_1 \right) -G\left( H_2 \right) \right)
   \cdot u \right]|| \leq ||\left( G\left( H_1 \right) -G\left( H_2 \right) \right) \cdot
   u||,
\end{equation*}
but for all $i=1,\ldots,n$:
\begin{equation}
||\left( G\left( H_1 \right) -G\left( H_2 \right) \right) \cdot
u||\leq ||G_i\left( H_1 \right) -G_i\left( H_2 \right)|| \, ||u||,
\label{eq:firststep}
\end{equation}
and by Proposition~\ref{powerseries} $||G_i\left( H_1 \right) -G_i\left( H_2 \right)||\leq \frac{||G_i||_{r_i}}{r_i}||H_1-H_2||$. Then setting $\mu=||u|| \frac{\max M_i}{\min r_i}$, from~\eqref{eq:firststep}
 we conclude that
\begin{equation*}
||H_1-H_2||\leq \mu ||H_1-H_2||.
\end{equation*}
By hypothesis $||u|| <\frac r M$, then $\mu <1$, from which we conclude that $H_1=H_2$.

\indent We now prove existence. Since $G_i$ are
convergent, $D^{\alpha}G_i$ also are convergent  and
Proposition~\ref{Taylorformulas} gives the following estimate
\begin{equation}
||D^{\alpha}G_i||_{r_i}\leq \frac{M_i}{r_i^{\lvert \alpha
\rvert}}, \quad \forall \alpha \in \mathbb{N}^n \label{eq:estimcauchy} \, ,
\end{equation}
which together with the non--expanding
 property of $\Lambda$
allows us to prove that for all $N \geq 1$ and all $\vartheta \in
\mathcal{T}_{N}$:
\begin{equation*}
 \Big\lvert \Big\lvert \vall{\vartheta}{u,w}{\Lambda} \Big\rvert \Big\rvert \leq
  ||u||^{N}\frac{M^N}{r^{N-1}} \, .
\end{equation*}
Let us define $H^{(0)}\left( u,w \right)=\Lambda w$ and
$H^{(j)}\left( u,w \right)=\Lambda w+\sum_{N=1}^j \sum_{\vartheta
\in \mathcal{T}_{N}}\vall{\vartheta}{u,w}{\Lambda}$, clearly
$H^{(j)}\rightarrow H$ as $j \rightarrow \infty$ and it is easy to
check that:
\begin{equation*}
\Big\lvert \Big\lvert H^{(j)}-\Lambda \left[ w-G\left( H^{(j)}
\right) \cdot u \right]\Big\rvert\Big\rvert \leq r
\left(||u||\frac{M}{r}\right)^{j+1}
\end{equation*}
which tends to $0$ as $j \rightarrow \infty$.

\endproof

We give now the proof of Corollary~\ref{existence}. This one follows
closely the one of Theorem~\ref{nonanex} in particular the
uniqueness statement, so we will outline only the main differences
w.r.t to the previous proof.

\proof The hypothesis on $G^{\prime}_i$ gives an estimate similar
to~\eqref{eq:estimcauchy}, then by induction on $N$ it is easy to
prove that for all $N\geq 1$ and $\vartheta \in \mathcal{T}_{N}$
one has
\begin{equation*}
\Big\lvert\Big\lvert\val{\vartheta}{u^{\prime},w^{\prime}}\Big\rvert\Big\rvert\leq ||u^{\prime}||^{N}\frac{M^N}{r^{N-1}}.
\end{equation*}
Then if $||u^{\prime}||<\frac{r}{M}$ series~\eqref{solution} converges and if
$||w^{\prime}||\leq r$, $H\left(u^{\prime},w^{\prime} \right) \in B_0 \left( 0, r
\right)$, in fact
\begin{equation*}
||H||=\Big\lvert\Big\lvert w^{\prime}+\sum_N \sum_{\vartheta \in
\mathcal{T}_{N}} \val{\vartheta}{u^{\prime},w^{\prime}}\Big\rvert\Big\rvert \leq
\sup \Big \{ ||w^{\prime}||,\sup_{N,\vartheta \in \mathcal{T}_{N}}
\Big\lvert\Big\lvert \val{\vartheta}{u^{\prime},w^{\prime}}\Big\rvert\Big\rvert \Big
\} < r.
\end{equation*}
Now introducing $H^{(0)}\left( u^{\prime},w^{\prime} \right)=w$ and
$H^{(j)}\left( u^{\prime},w^{\prime} \right)=w^{\prime}+\sum_{N=1}^j \sum_{\vartheta \in
\mathcal{T}_{N}}\val{\vartheta}{u^{\prime},w^{\prime}}$, one can easily prove that
$H^{(j)}\rightarrow H$ as $j \rightarrow \infty$.
\endproof

\begin{remark}
In the simplest case $n=l=1$, namely $u,w\in k$ and  $G\in
k\left[\left[ X \right]\right]$, the solution given
by~\eqref{solution} coincides with the classical one of
Lagrange~\eqref{eq:lagrangeorig}. One can prove this fact either
using the uniqueness of the Taylor development or by direct
calculation showing that for all positive integer $N\geq 1$ we have
\begin{equation*}
\sum_{\vartheta \in \mathcal{T}_{N}}\val{\vartheta}{u,w}
=\frac{u^N}{N!}\frac{d^{N-1}}{dw^{N-1}}\left[ G(w)\right]^N .
\end{equation*}
In the other cases formula~\eqref{solution} is the natural
generalization of~\eqref{eq:lagrangeorig}.

\end{remark}

\begin{remark}
\label{rem:analytic} Series~\eqref{solution} is an analytic 
function of $u,w$, but it is not explicitly written as $u$--power series. 
We claim that introducing {\em labeled rooted trees} we can rewrite~\eqref{solution} explicitly as a $u$--power series.

\end{remark}

\section{The Non--analytic Siegel center problem.}
\label{sec:someappl}

In this part we show that the problem of the conjugation of a (formal) 
germ of a given function with its linear part near a fixed point 
(the so called {\em Siegel Center Problem}) can be solved applying
 Theorem~\ref{nonanex} to the field of (formal) power series. The Siegel
 Center Problem is a particular case of first order semilinear $q$--difference 
equation, but our results apply to general first order semilinear $q$--difference and differential equations (see next section).

\subsection{Notations and Statement of the Problem}
\label{subsec:siegeln}

Let $\alpha =(\alpha_1,\dots,\alpha_n)\in \N^n$, $\lambda=(\lambda_1,\dots,\lambda_n)\in \C^n$ with $\lambda_i \neq \lambda_j$ if $i\neq j$, we will use the compact notation $\lambda^{\alpha}=\lambda_1^{\alpha_1}\dots \lambda_n^{\alpha_n}$ and $|\alpha|=\alpha_1+\dots+\alpha_n$; we will
 denote the diagonal $n\times n$ matrix with $\lambda_i$ at the $(i,i)$--th
 place with $\diag{\lambda}$.

Let $V=\C^n\left[\left[z_1,\dots,z_n\right]\right]$ be the vector space of the 
formal power series in the $n$ 
variables $z_1,\dots,z_n$ with coefficients in $\mathbb{C}^n$: ${\bf f}\in V$ if
\begin{equation*}
{\bf f}=\sum_{\alpha \in \mathbb{N}^n} {\bf f}_{\alpha}z^{\alpha} ,\quad {\bf f}_{\alpha}\in \mathbb{C}^n \text{ and } z^{\alpha}=z_1^{\alpha_1}\dots z_n^{\alpha_n} \; \forall \; \alpha =(\alpha_1,\dots,\alpha_n)\in \mathbb{N}^n \, .
\end{equation*}

We consider $V$ endowed with the ultrametric absolute norm induced by 
the ${\bf z}$--adic valuation (${\bf z}=(z_1,\dots,z_n)$):
$||{\bf f}||=2^{-v({\bf f})}$, where 
$v({\bf f})=\inf \{|\alpha|, \alpha \in \mathbb{N}^n: {\bf f}_{\alpha}\neq 0\}$, and for
 any positive integer $j$ we denote by $V_j = \{{\bf f}\in V :v({\bf f})> j\}$.

Let $\mathcal{C}$ be a Class (that we will define later, see
 paragraph~\ref{classiMn}) of formal power series, closed w.r.t. the 
(formal) derivation, the composition and where, roughly speaking, the
 (formal) Taylor series makes sense. One can think for example to 
the Class of germs of analytic diffeomorphisms of $(\mathbb{C}^n,0)$ or
 Gevrey--$s$ Classes,
 in fact we will see that our classes will contain these special cases. 

Let $A\in GL(n,\mathbb{C})$ and assume $A$ to be diagonal~\footnote{The case 
$A$ non--diagonal need some special attentions, see~\cite{Herman} Proposition 
3 page 143 and~\cite{Yoccoz} Appendix 1.} with all the eigenvalues distinct. Let $\mathcal{C}_1$ and 
$\mathcal{C}_2$ be two classes as stated before, then the {\em Siegel center problem}
 can be formulated as follow~\cite{Herman,CarlettiMarmi}:
\begin{quotation}
Let ${\bf F}({\bf z})=A{\bf z}+{\bf f}({\bf z})\in \mathcal{C}_1$, ${\bf f}\in \mathcal{C}_1\cap V_1$, find necessary and sufficient conditions on $A$ to {\em linearize} in $\mathcal{C}_2$ ${\bf F}$, namely find ${\bf H}\in \mathcal{C}_2\cap V_0$ (the {\em linearization}) solving: 
\begin{equation}
\label{eq:siegelformal}
{\bf F}\circ {\bf H}({\bf z})={\bf H}(A{\bf z}) \, .
\end{equation}
\end{quotation}

Let $A=\diag{\lambda}$ we introduce the operator $D_{\lambda}:V\rightarrow V$:
\begin{equation}
\label{eq:defoperator}
D_{\lambda}{\bf g}({\bf z})={\bf g}(A{\bf z})-A{\bf g}({\bf z}) \, ,
\end{equation}
for any ${\bf g}({\bf z})\in V$. We remark that the action of $D_{\lambda}$  on the monomial ${\bf v}z^{\alpha}$, for any ${\bf v}\in\C^n$ and any $\alpha \in\N^n$, is given by:
\begin{equation}
\label{eq:actionmonomial}
D_{\lambda}{\bf v}(z^{\alpha})=(\Omega_{\alpha}{\bf v})z^{\alpha} \, ,
\end{equation}
where the matrix $\Omega_{\alpha}$ is defined by $\Omega_{\alpha}= {\it diag}(\lambda^{\alpha}-\lambda_1,\dots,\lambda^{\alpha}-\lambda_n)$ and $(\Omega_{\alpha}{\bf v})=\sum_{i=1}^{n}(\lambda^{\alpha}-\lambda_i)v_i$ is the matrix--vector product.
 
Let $A=\diag{\lambda}\in GL(n,\mathbb{C})$, we say that $A$ is 
{\em resonant} if there exist $\alpha \in\N^n$, $|\alpha|\geq 2$
 and $j\in \{ 1,\dots,n\}$ such that:
\begin{equation}
\label{eq:resonant}
\lambda^{\alpha}-\lambda_j = 0 \, .
\end{equation}
If~\eqref{eq:resonant} doesn't hold, we say that $A$ is {\em non--resonant}.

If $A$ is non--resonant then $D_{\lambda}$ is invertible on $V_1$ (namely $|\alpha| \geq 2$), it is non--expanding ($||D_{\lambda}{\bf g}||\leq ||{\bf g}||$), and clearly $V$--additive ($D_{\lambda}({\bf f+g})=D_{\lambda}{\bf f}+D_{\lambda}{\bf g}$) and $\C$--linear. Then we claim that the Siegel center problem~\eqref{eq:siegelformal} is equivalent to solve the functional equation:
\begin{equation}
\label{eq:siegelsecvers}
D_{\lambda}{\bf h}={\bf f}\circ ( {\bf z}+{\bf h})
\end{equation}
where ${\bf z}$ is the identity formal power series, ${\bf f}\in \mathcal{C}_1\cap V_1$ and ${\bf h}\in \mathcal{C}_2\cap V_1$. In fact from~\eqref{eq:siegelformal} we see that the linear part of ${\bf H}$ doesn't play any role, so we can choose ${\bf H}$ tangent to the identity (this normalization assures the uniqueness of the linearization): ${\bf H}({\bf z})={\bf z}+{\bf h}({\bf z})$, with ${\bf h}\in \mathcal{C}_2\cap V_1$. But then~\eqref{eq:siegelformal} can be rewritten as:
\begin{equation}
\label{eq:siegelstep1}
{\bf h}(A{\bf z})-A{\bf h}({\bf z})={\bf f}\circ ({\bf z} + {\bf h}({\bf z})) \, ,
\end{equation}
and replacing the left hand side with the operator $D_{\lambda}{\bf h}$ we obtain~\eqref{eq:siegelsecvers}.

\indent Given ${\bf f}\in V_1$ we consider the function $G_{\bf f}({\bf g})={\bf f}\circ ({\bf z} +{\bf g})$, for all ${\bf g}\in V_1$, assume that we can invert the operator $D_{\lambda}$ (because $v({\bf f}({\bf z} +{\bf g}))\geq v({\bf f})$ we can invert $D_{\lambda}$ whenever ${\bf f}\in V_1$), we can rewrite~\eqref{eq:siegelsecvers} as:
\begin{equation}
\label{eq:siegeln}
{\bf h} = D_{\lambda}^{-1} \left(G_{\bf f}({\bf h})\right) \, , \end{equation}
which is a particular case of the non--analytic multidimensional Lagrange inversion formula~\eqref{multilagrange} with $u=1$, $w=0$ and $\Lambda=D^{-1}_{\lambda}$. The following Lemma assures that $G_{\bf f}$ verifies the hypotheses of Theorem~\ref{nonanex}.

\begin{lemma}
\label{lem:gfv}Given ${\bf f}\in V_1$ the composition ${\bf f}\circ ({\bf z} +{\bf v})$ defines a power series
$G_{\bf f}\in V_1\left[\left[ v_1,\dots , v_n\right]\right]$ convergent on $B(0,1/2)=V_0$ and
\begin{equation}
\label{eq:lemmahypot}
G_{\bf f}({\bf v})=\sum_{\beta \in\mathbb{N}^n}{\bf g}_{\beta}({\bf f})v^{\beta}
\end{equation}
where (if we define ${\bf f}_{\beta}=0$ for $|\beta|=0$ and $|\beta|=1$) series ${\bf g}_{\beta}({\bf f})\in V$ are given by
\begin{equation}
\label{eq:defofg}
{\bf g}_{\beta}({\bf f})({\bf z})=\sum_{\alpha \in\mathbb{N}^n}{\bf f}_{\alpha +\beta}
\binom{\alpha  +\beta}{\beta}z^{\alpha}
\end{equation}
Here we used the compact notations $\alpha !=\alpha_1 !\dots  \alpha_n!$ and
 $\binom{\alpha +\beta}{\beta}=\frac{(\alpha + \beta)!}{\alpha !\beta !}$, for 
$\alpha=(\alpha _1,\dots,\alpha_n)\in \mathbb{N}^n$ and $\beta \in \mathbb{N}^n$.
Moreover one has
\begin{equation}
\label{eq:otherforg}
||{\bf g}_{0}||=||{\bf f}||, \, ||{\bf g}_{\beta}||=2||{\bf f}|| \text{ for any $|\beta|=1$ and } \, ||{\bf g}_{\beta}||\leq 4||{\bf f}|| \text{ for any $|\beta|\geq 2$};
\end{equation}
thus
\begin{equation}
\label{eq:lastforG}
||G_{\bf f}||_{1/2}\leq ||{\bf f}||.
\end{equation}
\end{lemma}

The proof is straightforward and we omit it.

We can thus apply Theorem~\ref{nonanex} to solve~\eqref{eq:siegeln} with
$u=1$, $w=0$, $G=G_{\bf f}$, $\Lambda =D^{-1}_{\lambda}$, $r=1/2$ and $M=1/4$, and the unique solution of~\eqref{eq:siegeln} is given by
\begin{equation}
\label{eq:solution}
{\bf h}=\sum_{N\geq 1}\sum_{\vartheta \in \forest{N}}\vall{\vartheta}{1,0}{D^{-1}_{\lambda}}.
\end{equation}

An explicit expression for the power series coefficients of ${\bf h}$ can be
obtained introducing {\em labeled rooted trees} (see Remark~\ref{rem:analytic}). Let us now explain how to do this. Let $N\geq 1$ and let $\vartheta$~\label{pg:label} be a rooted tree of order $N$; to any $v\in \vartheta$ we associate a {\em (node)--label}  $\alpha_v \in \N^n$ s.t. $|\alpha_v|\geq 2$ and to the line $\ell_v$ (exiting w.r.t. the partial order from $v$) we associate a {\em (line)--label} $\beta^{\ell_v}\in \N^n$ s.t. $|\beta^{\ell_v}|=1$. We define
 $\beta_v =\sum_{\ell \in L_v}\beta^{\ell}$ (so $\beta_v \in \N^n$ and 
$|\beta_{v}|=m_v$) and the {\em momentum} flowing through the line $\ell_v$: $\nu_{\ell_v}=\sum_{w\in \vartheta :w\leq v}(\alpha_w -\beta_w)$. When $v$
 will be root of the tree, we will also use the symbol $\nu_{\vartheta}$ 
(the {\em total momentum} of the tree) instead of $\nu_{\ell_{v}}$. It is 
trivial to show that the momentum function is increasing (w.r.t. the partial
 order of the tree), namely if $v$ is the root of a rooted labeled tree and
 if $v_i$ is any of the immediate predecessor of $v$, $v > v_i$, we have:
 $|\nu_v| > |\nu_{v_i}|$). Recalling that the order of $\vartheta$ is $N$ we have: $|\nu_{\vartheta}|\geq N+1$.

Let $N\geq 1$, $j\in \{1,\dots,n \}$ and $\alpha \in \N^n$, we finally define $\mathcal{T}_{N,\alpha ,j}$ to be the {\em forest of rooted labeled trees} of order $N$ with total momentum $\nu_{\vartheta}=\alpha$ and $\beta^{\ell_{v_1}}=e_j$ 
(being $e_j$ the vector with all zero entries but the $j$--th which is set 
equal to $1$) for the root line $\ell_{v_1}$.

We are now able to prove the following

\begin{proposition}
\label{thecoefficients}
Equation~\eqref{eq:siegeln} admits a unique solution  ${\bf h}\in V_1$,  
${\bf h}=\sum_{\alpha \in \N ^n}{\bf h}_{\alpha}z^{\alpha}$.
 For $|\alpha |\ge 2$ the $j$--th component of the coefficient 
${\bf h}_{\alpha}$ is given by~\footnote{Compare this expression 
with equation (3.7) of Proposition  3.1 in~\cite{ChierchiaFalcolini}.}:
\begin{equation}
\label{coefficients}
h_{\alpha ,j }=\sum_{N=1}^{|\alpha |-1}\sum_{\vartheta\in
\mathcal{T}_{N,\alpha ,j}}
((\Omega^{-1}_{\nu_{\ell_{v_1}}}
{\bf f}_{\alpha_{v_1}})\cdot\beta^{\ell_{v_1}}) \prod_{v\in \vartheta } \binom{\alpha_v}{\beta_v}
\prod_{\ell_w \in L_v}((\Omega^{-1}_{\nu_{\ell_w}}
{\bf f}_{\alpha_w})\cdot\beta^{\ell_w})\, ,
\end{equation}
where the last product has to be set equal to $1$ whenever $v$ is an end node ($L_v=\emptyset$).
\end{proposition}

\begin{remark}
By definition $\beta^{\ell_w}$, for any $w \in\vartheta$, has length $1$, so 
it coincides with an element of the canonical base. Then for $w\in\vartheta$ 
and any choice of the labels, such that $\beta^{\ell_w}=e_i$, the term 
$((\Omega^{-1}_{\nu_{\ell_w}}
{\bf f}_{\alpha_w})\cdot\beta^{\ell_w})$ is nothing else that
\begin{equation*}
((\Omega^{-1}_{\nu_{\ell_w}}
{\bf f}_{\alpha_w})\cdot \beta^{\ell_w})=\frac{1}{\lambda^{\nu_{\ell_w}}-\lambda_i}f_{\alpha_w,i} \, .
\end{equation*}
\end{remark}
\begin{remark}\label{rem:negativecomp}Even if all the nodes labels 
$\alpha_v$ have non--negative components, the momenta can have (several) 
negative components. More precisely, let $\vartheta \in \forest{N}$, if
 the order of the tree if big enough w.r.t. the dimension $n$ ($N\geq n$) 
then $\nu_{\vartheta}$ can have $n-1$ negative components and their sum 
can be equal to $1-N$, but we always have $|\nu_{\vartheta}|\geq N+1$. 
This fact reminds the definitions of the sets ${\bf N}_i$, ${\bf N}^{(m)}$ 
and ${\bf N}^{(m)}_+$ of~\cite{Bruno}.\end{remark}

\proof 
Let $\vartheta \in \mathcal{T}_{N}$, for $N\ge 1$ and $\alpha \in \N^n$ 
such that $|\alpha|\geq N+1$. Let us define ${\it Val}(\vartheta)=
((\Omega^{-1}_{\nu_{\ell_{v_1}}}
{\bf f}_{\alpha_{v_1}})\cdot \beta^{\ell_{v_1}}) \prod_{v\in \vartheta } \binom{\alpha_v}{\beta_v}
\prod_{\ell_w \in L_v}((\Omega^{-1}_{\nu_{\ell_w}}
{\bf f}_{\alpha_w})\cdot \beta^{\ell_w})$ and
\begin{equation}
\label{eq:firsteqa}
h^{\vartheta}_{|\alpha|,N,j}=\sum_{|\nu_{\vartheta}|=|\alpha|}{\it Val}(\vartheta) \, ,
\end{equation}
namely for a fixed tree, sum over all possible labels $\alpha_{v_i}$ and $\beta^{\ell_{v_i}}$, with $v_i$ in the tree, in such a way that the 
total momentum is fixed to $\alpha$ and the root line has label 
$\beta^{\ell_{v_1}}=e_j$. It is clear that 
\begin{equation}
\label{eq:seceqa}
{\bf h}=\sum_{|\alpha |\geq 2}z^{\alpha}\sum_{N=1}^{|\alpha|-1}\sum_{\vartheta \in \forest{N}}{\bf h}^{\vartheta}_{|\alpha|,N} \, ,
\end{equation}
being ${\bf h}^{\vartheta}_{|\alpha|,N}$ the vector whose $j$--th component is  $h^{\vartheta}_{|\alpha|,N,j}$. Let ${\bf h}^{\vartheta}_{N+1}=\sum_{|\alpha |\geq N+1}z^{\alpha}{\bf h}^{\vartheta}_{|\alpha|,N}$, clearly $||{\bf h}^{\vartheta}_{N+1}||\leq 2^{-(N+1)}$ and
\begin{equation}
\label{eq:tereqa}
{\bf h}=\sum_{N\geq 1}\sum_{\vartheta \in \forest{N}}{\bf h}^{\vartheta}_{N+1} \, .
\end{equation}
Actually in~\eqref{eq:seceqa} ${\bf h}$ is ordered with increasing powers
 of ${\bf z}$ whereas in~\eqref{eq:tereqa} with increasing order of trees.
 Convergence in $V_1$ for~\eqref{eq:tereqa} is assured from the estimate
 $||{\bf h}^{\vartheta}_{N+1}||\leq 2^{-(N+1)}$ and uniform convergence 
assures that~\eqref{eq:seceqa} and~\eqref{eq:tereqa} coincide.

We claim that by induction on the order of the tree, we can prove that for
 all $N\geq 1$ and all $\vartheta \in \forest{N}$ :
\begin{equation}
\label{eq:induct}
{\bf h}^{\vartheta}_{N+1}=\vall{\vartheta}{1,0}{D^{-1}_{\lambda}} \, ,
\end{equation}
where $\vall{\vartheta}{1,0}{D^{-1}_{\lambda}}$ has been defined 
in~\eqref{defvallambda}, taking $\Lambda =D^{-1}_{\lambda}$ and $G=G_{\bf f}$, thus establishing the equivalence of~\eqref{coefficients} and~\eqref{eq:solution}.
\endproof

\subsection{Some known results}
\label{subsec:soemknownres}

If both Classes $\mathcal{C}_1$ and $\mathcal{C}_2$ are $V_1$ 
(which verify the hypotheses of stability w.r.t. the derivation, 
closeness w.r.t. the composition, and the formal Taylor series 
makes sense), then the {\em Formal} Siegel Center Problem has a 
solution if the linear part of ${\bf F}$ is non--resonant.

A matrix $A=\diag{\lambda}\in GL(n,\mathbb{C})$, is in the 
{\em Poincar\'e domain} if
\begin{equation}
\label{eq:poincaredomain}
\sup_{1\leq j\leq n}|\lambda_j|<1 \quad \text{or} \quad \sup_{1\leq j\leq n}|\lambda_j^{-1}|>1 \, ,
\end{equation}
if $A$ doesn't belong to the Poincar\'e domain it will be in 
the {\em Siegel domain}.

In the {\em Analytic} case (both $\mathcal{C}_1$ and $\mathcal{C}_2$ 
are the ring of the germs of analytic diffeomorphisms of $(\mathbb{C}^n,0)$)
, let $A$ be the derivative of ${\bf F}$ at the origin, then if $A$
 is non--resonant and it is in the {\em Poincar\'e domain}, the 
Analytic Siegel Center Problem has a solution~\cite{Poincare1,Koenigs}
 (see also~\cite{Herman} and references therein). Moreover if $A$
 is resonant and in the Poincar\'e domain, but ${\bf F}$ is formally 
linearizable, then ${\bf F}$ is analytically linearizable.

If $A$ is in the Siegel domain the problem is harder, but we can
 nevertheless have a solution of the Analytic Siegel Center Problem, 
introducing some new condition on $A$. Let $p\in \N$, $p\geq 2$, and let us define
\begin{equation}
\label{eq:piccolidivisori}
\Tilde\Omega (p)=\min_{1\leq j \leq n}\inf_{\alpha \in \Z^n: |\alpha|< p} |\lambda^{\alpha}-\lambda_j| \, ,
\end{equation}
we remark that even if $A$ is non--resonant, but in the Siegel domain, one 
has $\lim_{p\rightarrow \infty}\Tilde\Omega (p)=0$; this is the
 so called {\em small divisors problem}, the main obstruction to 
the solution of equation~\eqref{eq:siegelformal}. A non--resonant 
matrix $A$ verifies a {\em Bruno condition}~\footnote{See~\cite{Bruno,Russmann}
 for various equivalent formulations of this condition.} if there
 exists an increasing sequence of natural numbers $(p_k)_k$ such that
\begin{equation}
\label{eq:brunocondition}
\sum_{k=0}^{+\infty} \frac{\log \Tilde{\Omega}^{-1}(p_{k+1})}{p_{k}}<+\infty \, .
\end{equation}
Then if $A$ satisfies a Bruno Condition, the germ is analytically linearizable~\cite{Bruno,Russmann}.

For the $1$--dimensional Analytic Siegel Center Problem, 
Yoccoz~\cite{Yoccoz} proved that the Bruno condition is 
necessary and sufficient to linearize analytically any univalent germs
with fixed linear part, in this case the Bruno 
condition reduces to the convergence of the series
\begin{equation}
\label{eq:1dimbruno}
\sum_{k=0}^{+\infty} \frac{\log q_{k+1}}{q_{k}}<+\infty \, ,
\end{equation}
where $(q_k)_k$ is the sequence of the convergent to 
$\omega \in \R \setminus \Q$ such that $\lambda=e^{2\pi i\omega }$.

\subsection{A new result: ultradifferentiable Classes}
\label{classiMn}

Let $(M_{k})_{k\ge 1}$ be a sequence of positive real numbers such that: 
\label{pg:ultradiff}
\item{0)} $\inf_{k\ge 1} M_{k}^{1/k}>0$; 
\item{1)} There exists $C_{1}>0$ such that $M_{k+1}\le C_{1}^{k+1}M_{k}$ for all $k\ge 1$; 
\item{2)} The sequence $(M_{k})_{k\ge 1}$ is logarithmically convex;
\item{3)} $M_{k}M_{l}\le M_{k+l-1}$ for all $k,l\ge 1$.

We define the class $\class{M_k}\subset  \C^n\left[\left[z_1,\dots,z_n\right]\right]$ as the set of formal power series ${\bf f}=\sum {\bf f}_{\alpha}z^{\alpha}$ such that there exist $A,B$ positive constant, such that:
\begin{equation}
|{\bf f}_{\alpha}|\leq A B^{|\alpha|}M_{|\alpha|} \quad \forall \alpha \in\mathbb{N}^n \, .
\label{eq:classMk}
\end{equation}
The hypotheses on the sequence $(M_k)_k$ assure that $\class{M_k}$ 
is stable w.r.t. the (formal) derivation, w.r.t. the composition of 
formal power series and for every tensor built with element of the class, its contraction\footnote{This assures that any term of the Taylor series 
is well defined.} gives again an element of the class. For example 
if ${\bf f},{\bf g}\in \class{M_k}$ then also $d{\bf f}({\bf z})({\bf g}({\bf z}))$ belongs to the same class. 
\begin{remark}
Our classes include the Class of Gevrey--$s$ power series as a
 special case: $M_k=(k!)^s$. Also the ring of convergent (analytic)
 power series are trivially included.
\end{remark}

In~\cite{CarlettiMarmi} a similar problem was studied in the
 $1$--dimensional case. Here we will extend the results contained there to the case of dimension
 $n \geq 1$. The main result will be the following Theorem

\begin{theorem}
\label{maintheorem}
Let $(\lambda_1,\dots,\lambda_n)\in \C^n$, $|\lambda_i|=1$ for $i=1,\dots ,n$, and let $A=\diag{\lambda}$ be non--resonant,
 $(M_k)_k$ and $(N_k)_k$ be sequences verifying hypotheses 0)--3). 
Let ${\bf F}\in V_0$, s.t. ${\bf F}({\bf z})=A{\bf z}+{\bf f}({\bf z})$ 
where ${\bf f}\in V_1$. Then
\begin{enumerate}
\item If moreover ${\bf F}\in \class{M_k}\cap V_0$ and $A$ verifies a
 Bruno condition~\eqref{eq:brunocondition}, then also the linearization 
${\bf H}$ belongs to $\class{M_k}\cap V_0$.

\item If ${\bf F}$ is a germ of analytic diffeomorphisms of $(\C^n,0)$
and there exists an increasing sequence of integer numbers $(p_k)_k$
 such that $A$ verifies:
\begin{equation}
\label{eq:firstnewcondition}
\limsup_{|\alpha |\rightarrow +\infty} \left(  2\sum_{m=0}^{\kappa (\alpha)} \frac{\log \Omega^{-1}(p_{m+1})}{p_m}-\frac{1}{|\alpha|}\log N_{|\alpha|} \right) <+\infty \, ,
\end{equation}
where $\kappa(\alpha)$ is the integer defined by: $p_{\kappa(\alpha)}\leq |\alpha| <p_{\kappa(\alpha)+1}$, then the linearization ${\bf H}$ belongs 
to $\class{N_k}\cap V_0$.

\item If ${\bf F}\in \class{M_k}\cap V_0$, the sequence $(M_k)_k$ 
is asymptotically bounded by the sequence $(N_k)_k$ (namely $N_k \geq M_k$ 
for all sufficiently large $k$) and there exists an increasing sequence 
of integer numbers $(p_k)_k$ such that $A$ verifies:
\begin{equation}
\label{eq:secondnewcondition}
\limsup_{|\alpha |\rightarrow +\infty} \left(  2\sum_{m=0}^{\kappa (\alpha)} \frac{\log \Omega^{-1}(p_{m+1})}{p_m}-\frac{1}{|\alpha|}\log \frac{N_{|\alpha|}}{M_{|\alpha|}} \right) <+\infty \, ,
\end{equation}
where $\kappa(\alpha)$ is the integer defined by: $p_{\kappa(\alpha)}\leq |\alpha| <p_{\kappa(\alpha)+1}$, then the linearization ${\bf H}$ belongs to $\class{N_k}\cap V_0$.

\end{enumerate}
\end{theorem}
The proof of Theorem~\ref{maintheorem} will be done in 
section~\ref{subsub:proofoftheorem}, before we make some 
remarks and we prove some preliminary lemmata.
\begin{remark}
The new arithmetical conditions~\eqref{eq:firstnewcondition}
 and~\eqref{eq:secondnewcondition} are generally weaker than the 
Bruno condition~\eqref{eq:brunocondition}. Theorem~\ref{maintheorem} 
is the natural generalization of results of~\cite{CarlettiMarmi} 
(compare it with Theorem 2.3 of~\cite{CarlettiMarmi}): condition~\eqref{eq:secondnewcondition} (respectively~\eqref{eq:firstnewcondition}) reduces to condition (2.10) (respectively condition (2.9)) of~\cite{CarlettiMarmi} except for the factor in front of the sum: here we have $2$ instead of $1$ in~\cite{CarlettiMarmi}. This is due to the better control of small denominators one can achieve using continued fractions and Davie's counting lemma~\cite{Davie} as explained in~\cite{CarlettiMarmi}.  
\end{remark}

\subsubsection{Some preliminaries}
\label{subsub:prelim}
Let $\omega_j \in (-1/2,1/2)\setminus \Q$ for $j=1,\dots,n$ and assume $\lambda=(\lambda_1,\dots,\lambda_n)$ with $\lambda_j=e^{2\pi i \omega_j}$. Let $\nu_{\ell}$, be the momentum of a line of a rooted labeled tree which contributes with a small denominator of the form: $|\lambda^{\nu_{\ell}}-\lambda_j|=2|\sin \pi(\nu_{\ell}\cdot  \omega - \omega_j)|$, where $\nu_{\ell} \cdot \omega=\sum_{j=1}^{n}{\nu_{\ell}}_j \omega_j$. Then using
\begin{equation*}
2|x|\leq|\sin \pi x|\leq \pi |x| \quad \forall |x|\leq \frac 1 2 \, ,
\end{equation*}
we claim that the contribution of the small denominator is equivalent 
to $\{\nu_{\ell}\cdot  \omega - \omega_j\}$, where $\{x\}$ denotes the 
distance of $x$ from its nearest integer. 

Let $p\in\N$, $p \geq 2$, and $\Omega(p)=\min \{ \{\nu \cdot \omega \}, \nu\in\Z^n : 0< |\nu|\leq p \}$. Let $(p_k)_k$ be an increasing sequence of positive integer and define (see~\cite{Bruno})
\begin{equation}
\label{eq:definphi}
\Phi^{(k)}(\nu)=\begin{cases}
                1 \quad \text{if } \{ \nu \cdot \omega \} < \frac 1 2 \Omega(p_k) \\
                0 \quad \text{if } \{ \nu \cdot \omega \} \geq \frac 1 2 \Omega(p_k) \, ,
               \end{cases}
\end{equation}
for any $\nu\in\Z^n \setminus 0$. By definition we trivially have $\Phi^{(k)}(\nu)=0$ for all $0\leq |\nu|\leq p_k$.

We define the following Bruno condition:
\begin{equation}
\label{eq:bruno2}
\sum_{k=0}^{+\infty} \frac{\log \Omega^{-1}(p_{k+1})}{p_{k}}<+\infty \, .
\end{equation}
and we claim that it is equivalent to~\eqref{eq:brunocondition}. We can prove the following
\begin{lemma}
\label{lemma:brj1}
Let $\nu_1 \in \Z^n$ such that $\Phi^{(k)}(\nu_1)=1$, for some $k$. Then for all $\nu_2 \in \Z^n$, such that $0 < | \nu_2 | \leq p_k$, we have $\Phi^{(k)}(\nu_1-\nu_2)=0$.
\end{lemma}
The proof follows closely the one of Lemma 10 p.218 of~\cite{Bruno} and we don't prove it.

Let $\ell$ be a line of a rooted labeled tree, $\vartheta$, let us introduce the notion of {\em scale} of the line $\ell$. Let $\nu_{\ell}$ be the momentum of the line and let us define $\bar{\nu}_{\ell}=\nu_{\ell}-\beta^{\ell}$. Let $(p_k)_k$ be an increasing sequence of positive integer, then for any $k \geq 0$ we define:
\begin{equation}
\label{eq:scaleline}
s_{\ell}(k)=\begin{cases}
             1 \quad \text{if } \frac 1 2 \Omega(p_{k+1})\leq \{ \bar{\nu}_{\ell}\cdot \omega \}< \frac 1 2 \Omega(p_{k}) \\
             0  \quad \text{otherwise . }
             \end{cases}
\end{equation}
For short we will say that a line $\ell$ is on scale $k$ if $s_{\ell}(k)=1$. Let us define $N_k(\vartheta)$ be the number of line on scale $1$ in the rooted labeled tree $\vartheta$. We can now prove the following Lemma, which roughly speaking says that the number of "bad" (too small) denominators is not too big, whereas Lemma~\ref{lemma:brj1} says that they do not occur so often.
\begin{lemma}[Bruno's Counting lemma]
\label{lemma:conteggio}
The number of lines on scale $k$ in a rooted labeled tree verifies the following bound:
\begin{equation}
\label{eq:boundcount}
N_k(\vartheta)\leq 
\begin{cases}
0 \quad &\text{if }|\bar{\nu}_{\vartheta}|< p_k \\
2\Big \lfloor \frac{|\bar{\nu}_{\vartheta}|}{p_k} \Big \rfloor -1\quad &\text{if }|\bar{\nu}_{\vartheta}|\geq p_k \, .
\end{cases}
\end{equation}
where $\lfloor x\rfloor$ denotes the integer part of the real number $x$. We recall that $\nu_{\vartheta}$ is the total momentum of the tree and $\bar{\nu}_{\vartheta}=\nu_{\vartheta}-\beta^{\ell_{v_1}}$, being $\ell_{v_1}$ the root line of $\vartheta$.
\end{lemma}

Our proof follows the original one of Bruno but exploiting the tree formalism, the interested reader can find this proof in appendix~\ref{sec:brunocl}.

\subsubsection{Proof of Theorem~\ref{maintheorem}}
\label{subsub:proofoftheorem}

We are now able to prove the main Theorem, we will prove only point 3
 which clearly implies point 1 (choosing $M_k=N_k$ for all $k$) and 
point 2 (choosing $M_k=C^k$ for all $k$ and some constant $C>0$).

Let us then assume that ${\bf F}\in V_0$, is of the form 
${\bf F}({\bf z})=A{\bf z}+{\bf f}({\bf z})$ where
 $\lambda =(\lambda_1, \dots, \lambda_n)\in \C^n$, $|\lambda_i|=1$ 
for $i=1,\dots ,n$, $A=\diag{\lambda}$ and ${\bf f}\in \class{M_k}\cap V_1$.

For a fixed rooted labeled tree of order $N \geq 1$ with total
 momentum equals to $\alpha \in\N^n$, $|\alpha |\geq 2$, we 
consider the following term of equation~\eqref{coefficients}:
\begin{equation}
\label{eq:coefffix}
((\Omega^{-1}_{\nu_{\ell_{v_1}}}
{\bf f}_{\alpha_{v_1}})\cdot \beta^{\ell_{v_1}}) \prod_{v\in \vartheta }
\prod_{\ell_w \in L_v}((\Omega^{-1}_{\nu_{\ell_w}}
{\bf f}_{\alpha_w})\cdot \beta^{\ell_w}) \, .
\end{equation}
Recalling the definition of scale and the definition of number of lines on scale $k$ we can bound~\eqref{eq:coefffix} with
\begin{equation*}
\Big\lvert ((\Omega^{-1}_{\nu_{\ell_{v_1}}}
{\bf f}_{\alpha_{v_1}})\cdot \beta^{\ell_{v_1}}) \prod_{v\in \vartheta }
\prod_{\ell_w \in L_v}((\Omega^{-1}_{\nu_{\ell_w}}
{\bf f}_{\alpha_w})\cdot \beta^{\ell_w})
\Big\rvert \leq \prod_{m=0}^{\kappa (\alpha)}\left[2 \Omega^{-1}(p_{m+1})\right]^{N_m(\vartheta)}\prod_{v\in \vartheta}|{\bf f}_{\alpha_v}| \, ,
\end{equation*}
where $\kappa (\alpha)$ is the integer defined by: $p_{\kappa (\alpha)}\leq |\alpha|<p_{\kappa (\alpha)+1}$.

Using hypothesis 3. of paragraph~\ref{classiMn} on the sequence $(M_k)_k$ and the hypothesis ${\bf f}\in \class{M_k}\cap V_1$, we will obtain for some positive constant $A,B$ the bound:
\begin{equation*}
\prod_{v\in \vartheta}|{\bf f}_{\alpha_v}|  \leq \prod_{v\in \vartheta}AB^{|\alpha_v|}M_{|\alpha_v|}\leq A^N B^{\sum_{v\in\vartheta}|\alpha_v|}M_{\sum_{v\in\vartheta}|\alpha_v|-(N-1)} \, ,
\end{equation*}
by definition of total momentum: $\sum_{v\in\vartheta}|\alpha_v|-(N-1)=|\nu_{\vartheta}|$, which has been fixed to $|\alpha|$, then
\begin{equation*}
\prod_{v\in \vartheta}|{\bf f}_{\alpha_v}| \leq B^{-1}(AB)^N B^{|\alpha|}M_{|\alpha|} \, .
\end{equation*}
Using finally the bound of the Counting Lemma~\ref{lemma:conteggio} we get:
\begin{equation}
\label{eq:lafinale}
\log \frac{|{\bf h}_{\alpha}|}{N_{|\alpha|}} \leq |\alpha|\log C+\log \frac{M_{|\alpha|}}{N_{|\alpha|}}+2|\alpha|\sum_{m=0}^{\kappa(\alpha)}\frac{\log \Omega^{-1}(p_{m+1})}{p_m} \, ,
\end{equation}
for some positive constant $C$. Dividing~\eqref{eq:lafinale}  by $|\alpha|$ and passing to the limit superior we get the thesis.

\begin{remark}
Let us consider a particular $1$--dimensional Siegel--Schr\"oder center problem with a germ of the form: $F_{(k)}(z)=\lambda z\left( 1-\frac{z^k}{k}\right)$ for some integer $k\geq 1$, $\lambda =e^{2 \pi i \omega}$, $\omega \in \R \setminus \Q$. Let us call $R_{(k)}(\omega)$ the radius of convergence of the unique linearization associated to $F_{(k)}$. Then an easy adaptation of Theorem~\ref{maintheorem} case 1) with $\mathcal{C}_{(M_k)}=z\C \{ z \}$, allows us to prove:
\begin{equation*}
 \log R_{(k)}(\omega) \geq -\frac{1}{k}B(k\omega)+\frac{\log k - C_k}{k}\, , 
\end{equation*}
for some constant $C_k$ (depending on $k$ but independent of $\omega$).

This can explain the $1/k$--periodicity of $R_{(k)}(\omega)$, as a function of $\omega$, showed in Figures $5$ and $7$ of~\cite{Marmi2}.
\end{remark}

\section{Linearization of non--analytic vector fields}
\label{sec:nonanvf}

The aim of this section is to extend the analytic results of Bruno about 
the linearization of an analytic vector field near a singular point, to the case of ultradifferentiable vector fields. 
We will show that this problem can be put in the framework of 
Theorem~\ref{nonanex} and then obtaining an explicit (i.e. non--recursive) 
expression for the change of variables (the {\em linearization}) 
in which the vector field has a simpler form.

Our aim is also to point out the strong similarities of this problem with
 the Siegel Center Problem, previously studied. In particular when both problems are put in the framework of the multidimensional non--analytic Lagrange inversion formula on the field of formal Laurent series, they give rise to (essentially) the same problem. For this reason most results will only be stated without proofs, these being very close to the proofs of the previous section.

\subsection{Notation and Statement of the Problem.}
\label{ssec:vfnotstat}

In this section we will use the same notations given at the beginning of section~\ref{classiMn}. Let $A\in GL(n,\mathbb{C})$ and assume $A$ to be diagonal. Let $\mathcal{C}_1$ and $\mathcal{C}_2$ be two classes of formal power series as defined before, then the {\em problem of the linearization of vector fields} can be formulated as follows:
\begin{quotation}
Let ${\bf F}({\bf z})=A{\bf z}+{\bf f}({\bf z})\in \mathcal{C}_1$, ${\bf f}\in \mathcal{C}_1\cap V_1$, and consider the following differential equation:
\begin{equation}
\label{eq:diffequa}
\dot {\bf z}=\frac{d {\bf z}}{dt}=A{\bf z}+{\bf f}({\bf z}) \, ,
\end{equation}
where $t$ denotes the time variable. Determine necessary and sufficient conditions on $A$ to find a change of variables in $\mathcal{C}_2\cap V_1$ (called the {\em linearization}) which leaves the singularity (${\bf z}=0$) fixed, doesn't change the linear part of ${\bf F}$ and allows to rewrite~\eqref{eq:diffequa} in a simpler form\footnote{Here we don't consider the most general case of looking for a change of coordinates which put~\eqref{eq:diffequa} in {\em normal form}, namely containing only resonant terms: $\dot w_i = w_i \sum_{\alpha \cdot \omega=0} g_{\alpha ,i} w^{\alpha}$. Our results will concern vector fields with non--resonant linear parts, so~\eqref{eq:equadiffformal} will be the normal form.}. Namely find ${\bf h}\in \mathcal{C}_2\cap V_1$, such that ${\bf z}= {\bf w}+{\bf h}({\bf w})$ and in the new variables ${\bf w}$, equation~\eqref{eq:diffequa} rewrites: 
\begin{equation}
\label{eq:equadiffformal}
\dot {\bf w}=A{\bf w}\, .
\end{equation}
\end{quotation}

Let $\omega=(\omega_1,\dots,\omega_n)\in\C^n$ and $A=\diag{\omega}$, we
 will say that $A$ is {\em resonant} if there exist $\alpha \in\Z^n$, with all positive component except at most one which can assume the value $-1$,
 $|\alpha |\geq 2$, and $j\in \{ 1,\dots,n\}$ such that: $\omega \cdot\alpha -\omega_j=0$, where $\omega \cdot\alpha=\sum_{i=1}^n \omega_i \alpha_i$ is the scalar product. Let $\alpha \in\N^n$ we introduce the diagonal matrix 
$\Omega^{\prime}_{\alpha}={\it diag}(\omega \cdot\alpha -\omega_1,\dots,\omega \cdot\alpha -\omega_n)$.

Let us introduce the operator\footnote{An equivalent definition will be: $D^{\prime}_{\omega}{\bf g}({\bf w})=\sum_{i=1}^n(A{\bf w})_i\partial_{w_i}{\bf g}({\bf w})-A{\bf g}({\bf w})=L_A \bf g$, the {\em Poisson bracket} of the linear field $A\bf w$ and $\bf g$.} $D^{\prime}_{\omega}:V\rightarrow V$ as follows
\begin{equation}
\label{eq:defoperatorvf}
D^{\prime}_{\omega}{\bf g}({\bf z})=\sum_{\alpha\in\N^n} (\Omega^{\prime}_{\alpha}{\bf g}_{\alpha})z^{\alpha}\, ,
\end{equation}
for any ${\bf g}({\bf z})\in V$, where $\Omega^{\prime}_{\alpha}\mathbf{g}_{\alpha}=\sum_{i=1}^n (\omega \cdot \alpha -\omega_i)g_{\alpha,i}$ is the matrix--vector product. 
If $A$ is non--resonant  then $D^{\prime}_{\omega}$ is invertible on $V_1$ (namely $|\alpha| \geq 2$), it is non--expanding and clearly $V$--additive and $\C$--linear. Then we claim that the linearization, ${\bf h}$, is solution of the functional equation:
\begin{equation}
\label{eq:vfsecvers}
D^{\prime}_{\omega}{\bf h}={\bf f}\circ ( z+{\bf h})\, ,
\end{equation}
where ${\bf f}\in \mathcal{C}_1\cap V_1$, ${\bf h}\in \mathcal{C}_2\cap V_1$ and $\mathbf{z}$ denotes the identity formal power series.

\indent 
Given ${\bf f}\in V_1$ we consider the function $G_{\bf f}({\bf h})={\bf f}\circ (z +{\bf h})$, assume $A$ to be non--resonant to invert the operator $D_{\omega}$, then we rewrite~\eqref{eq:vfsecvers} as:
\begin{equation}
\label{eq:vfn}
{\bf h} = {D^{\prime}}_{\omega}^{-1} \left(G_{\bf f}({\bf h})\right) \, , 
\end{equation}
which is a particular case of the non--analytic multidimensional Lagrange inversion formula~\eqref{multilagrange}. Apart from the different operator, this equation is the same of the Siegel Center Problem~\eqref{eq:siegelsecvers}. Lemma~\ref{lem:gfv} assures that $G_{\bf f}$ verifies the hypotheses of Theorem~\ref{nonanex} and thus we can apply it with: $u=1$, $w=0$, $G=G_{\bf f}$, $\Lambda={D^{\prime}}_{\omega}^{-1}$, $r=1/2$ and $M=1/4$. The unique solution of~\eqref{eq:vfn} is then given by:
\begin{equation}
\label{eq:solvfalberi}
{\bf h}=\sum_{N\geq 1}\sum_{\vartheta \in \forest{N}}\vall{\vartheta}{1,0}{{D^{\prime}}^{-1}_{\omega}}.
\end{equation}
Once again we can give an explicit expression for the linearization ${\bf h}$ using rooted labeled trees. Introducing the same labels as for the Siegel center problem (see page~\pageref{pg:label}) we can prove the following

\begin{proposition}
\label{thecoefficientsvf}
Equation~\eqref{eq:vfn} admits a unique solution  ${\bf h}\in V_1$,  
${\bf h}=\sum_{\alpha \in \N ^n}{\bf h}_{\alpha}z^{\alpha}$.
 For $|\alpha |\ge 2$ the $j$--th component of the coefficient ${\bf h}_{\alpha}$ is given by:
\begin{equation}
\label{coefficientsvf}
h_{\alpha ,j }=\sum_{N=1}^{|\alpha |-1}\sum_{\vartheta\in
\mathcal{T}_{N,\alpha ,j}}
(({\Omega^{\prime}}^{-1}_{\nu_{\ell_{v_1}}}
{\bf f}_{\alpha_{v_1}})\cdot \beta^{\ell_{v_1}}) \prod_{v\in \vartheta } \binom{\alpha_v}{\beta_v}
\prod_{\ell_w \in L_v}(({\Omega^{\prime}}^{-1}_{\nu_{\ell_w}}
{\bf f}_{\alpha_w})\cdot \beta^{\ell_w})\, ,
\end{equation}
where the last product has to be set equal to $1$ whenever $v$ is an end node ($L_v=\emptyset$).
\end{proposition}

We don't prove this Proposition (whose proof is the same as the one of Proposition~\ref{thecoefficients}); moreover we point out that because both problems give rise to (essentially) the same multidimensional non--analytic Lagrange inversion formula, we can pass from one solution to the other with very small changes: ${\Omega^{\prime}}^{-1}_{\alpha}$ instead of $\Omega^{-1}_{\alpha}$.

If both classes $\mathcal{C}_1$ and $\mathcal{C}_2$ are $V_1$ (formal case) then the linearization problem has solution if $A$ is non--resonant. In the analytic case we distinguish again the {\em Poincar\'e domain} (the {\em convex hull} of the $n$ complex points $\omega_1,\dots,\omega_n$ {\em doesn't contain the origin}) and the {\em Siegel domain} (if they are not in the Poincar\'e domain). In the first case, under a non--resonance condition, Poincar\'e proved that the vector field is analytically linearizable, and then Dulac, in the resonant case, proved the conjugation to a normal form. In the Siegel non--resonant case Bruno proved analytic linearizability~\cite{Bruno}, under the Bruno condition which reads:
\begin{equation}
\label{eq:brunovf}
\sum_{k \geq 0}\frac{\log \Hat{\Omega}^{-1}(p_{k+1})}{p_{k}}<+\infty \, ,
\end{equation}
for some increasing sequence of integer number $(p_k)_k$, where $\Hat{\Omega}(p)=\min\{ |\alpha \cdot \omega|: \alpha \in\N^n, \alpha \cdot \omega\neq 0, 0<|\alpha |<p, \text{at most one component $\alpha_i=-1$.}\}$.

We now extend this kind of results to the case of ultradifferentiable vector fields. Namely we consider two classes $\class{M_k}$ and $\class{N_k}$, defined as in section~\ref{classiMn} and we prove the following

\begin{theorem}
\label{the:mainvf}
Let $(\omega_1,\dots,\omega_n)\in \C^n$, $A=\diag{\omega}$ non--resonant,
 $(M_k)_k$ and $(N_k)_k$ be sequences verifying hypotheses 0)--3) of section~\ref{classiMn}. 
Let ${\bf F}\in V_0$, s.t. ${\bf F}({\bf z})=A{\bf z}+{\bf f}({\bf z})$ 
where ${\bf f}\in V_1$. Then
\begin{enumerate}
\item If moreover ${\bf F}\in \class{M_k}\cap V_0$ and $A$ verifies a
 Bruno condition~\eqref{eq:brunovf}, then also the linearization 
${\bf h}$ belongs to $\class{M_k}\cap V_1$.

\item If ${\bf F}$ is a germ of analytic diffeomorphisms of $(\C^n,0)$
and there exists an increasing sequence of integer numbers $(p_k)_k$
 such that $A$ verifies:
\begin{equation}
\label{eq:firstnewconditionvf}
\limsup_{|\alpha |\rightarrow +\infty} \left(  2\sum_{m=0}^{\kappa (\alpha)} \frac{\log \Hat{\Omega}^{-1}(p_{m+1})}{p_m}-\frac{1}{|\alpha|}\log N_{|\alpha|} \right) < +\infty\, ,
\end{equation}
where $\kappa(\alpha)$ is the integer defined by: $p_{\kappa(\alpha)}\leq |\alpha| <p_{\kappa(\alpha)+1}$, then the linearization ${\bf h}$ belongs 
to $\class{N_k}\cap V_1$.

\item If ${\bf F}\in \class{M_k}\cap V_0$, the sequence $(M_k)_k$ 
is asymptotically bounded by the sequence $(N_k)_k$ (namely $N_k \geq M_k$ 
for all sufficiently large $k$) and there exists an increasing sequence 
of integer numbers $(p_k)_k$ such that $A$ verifies:
\begin{equation}
\label{eq:secondnewconditionvf}
\limsup_{|\alpha |\rightarrow +\infty} \left(  2\sum_{m=0}^{\kappa (\alpha)} \frac{\log \Hat{\Omega}^{-1}(p_{m+1})}{p_m}-\frac{1}{|\alpha|}\log \frac{N_{|\alpha|}}{M_{|\alpha|}} \right) <+\infty\, ,
\end{equation}
where $\kappa(\alpha)$ is the integer defined by: $p_{\kappa(\alpha)}\leq |\alpha| <p_{\kappa(\alpha)+1}$, then the linearization ${\bf h}$ belongs to $\class{N_k}\cap V_1$.

\end{enumerate}
\end{theorem}

To prove this Theorem we will use again the majorant series method, the main step is to control the small denominators contributions. To do this, given an increasing sequence of positive integer $(p_k)_k$ we define a new {\em counting function}:
\begin{equation}
\label{eq:hatdefinphi}
\Hat{\Phi}^{(k)}(\nu)=\begin{cases}  
              1 \quad \text{if } | \nu \cdot \omega | < \frac 1 2 \Hat{\Omega}(p_k) \\
              0 \quad \text{if } | \nu \cdot \omega | \geq \frac 1 2 \Hat{\Omega}(p_k)\, .
               \end{cases}
\end{equation}
for any $\nu\in\Z^n\setminus 0$. By definition we trivially have $\Hat{\Phi}^{(k)}(\nu)=0$ for all $0< |\nu|\leq p_k$. Then we can prove the following Lemma (which will play the role of Lemma~\ref{lemma:brj1}):
\begin{lemma}
\label{lemma:brj1vf}
Let $\nu_1 \in \Z^n\setminus 0$ such that $\Hat{\Phi}^{(k)}(\nu_1)=1$, for some $k$. Then for all $\nu_2 \in \Z^n$, such that $0 < | \nu_2 | \leq p_k$, we have $\Hat{\Phi}^{(k)}(\nu_1-\nu_2)=0$.
\end{lemma}

We finally define a new notion of {\em scale}; let $\ell$ be a line of a rooted labeled tree, let $\nu_{\ell}$ be its momentum and let us recall that $\bar{\nu}_{\ell}=\nu_{\ell}-\beta^{\ell}$, then:
\begin{equation}
\label{eq:scalelinevf}
\Hat{s}_{\ell}(k)=\begin{cases}
             1 \quad \text{if } \frac 1 2 \Hat{\Omega}(p_{k+1})\leq | \bar{\nu}_{\ell}\cdot \omega |< \frac 1 2 \Hat{\Omega}(p_{k}) \\
             0  \quad \text{otherwise , }
             \end{cases}
\end{equation}
we will say that a line $\ell$ is on scale $k$ if $\Hat{s}_{\ell}(k)=1$. We can now prove the following Counting Lemma
\begin{lemma}[Bruno's Counting lemma $2^{\text{nd}}$ version]
\label{lemma:conteggiovf}
Let $\vartheta$ be a rooted labeled tree of order $N\geq 1$, let $k\geq 1$ be an integer and let $\Hat{N}_k(\vartheta)$ be the number of line on scale $k$ in the tree. Then the following bound holds:
\begin{equation}
\label{eq:boundcountvf}
\Hat{N}_k(\vartheta)\leq 
\begin{cases}
0 \quad &\text{if }|\bar{\nu}_{\vartheta}|< p_k \\
2\Big \lfloor \frac{|\bar{\nu}_{\vartheta}|}{p_k} \Big \rfloor -1\quad &\text{if }|\bar{\nu}_{\vartheta}|\geq p_k \, .
\end{cases}
\end{equation}
where $\lfloor x\rfloor$ denotes the integer part of the real number $x$. We recall that $\nu_{\vartheta}$ is the momentum of the root line and $\bar{\nu}_{\vartheta}=\nu_{\vartheta}-\beta^{\ell_{v_1}}$, being $\ell_{v_1}$ the root line.
\end{lemma}

We don't prove it because its proof is the same as the one of Lemma~\ref{lemma:conteggio} except for the use of Lemma~\ref{lemma:brj1vf} instead of Lemma~\ref{lemma:brj1}.

\subsection{Proof of Theorem~\ref{the:mainvf}}
\label{ssec:proofnv}
Once again we will prove only point 3 which clearly contains points 1 and 2 as special cases.

Assume ${\bf F}\in V_0$ of the form ${\bf F}({\bf z})=A{\bf z}+{\bf f}({\bf z})$ where $\omega =(\omega_1,\dots , \omega_n)\in\C^n$, $A=\diag{\omega}$ and ${\bf f}\in\class{M_k}\cap V_1$. Consider the contribution of a rooted labeled tree of order $N\geq 1$ with total momentum $\alpha\in\N^n$, $|\alpha|\geq 2$, given by~\eqref{coefficientsvf}:
\begin{equation}
\label{eq:equa1vf}
(({\Omega^{\prime}}^{-1}_{\nu_{\ell_{v_1}}}
{\bf f}_{\alpha_{v_1}})\cdot \beta^{\ell_{v_1}}) \prod_{v\in \vartheta } 
\prod_{\ell_w \in L_v}(({\Omega^{\prime}}^{-1}_{\nu_{\ell_w}}
{\bf f}_{\alpha_w})\cdot \beta^{\ell_w})\, .
\end{equation}
Follow closely the proof of Theorem~\ref{maintheorem} we can bound (use the definitions of $\Hat{N}_k(\vartheta)$ and of $\class{M_k}$) this contribution with:
\begin{equation}
\label{eq:eqavf2}
\Big\lvert (({\Omega^{\prime}}^{-1}_{\nu_{\ell_{v_1}}}
{\bf f}_{\alpha_{v_1}})\cdot \beta^{\ell_{v_1}}) \prod_{v\in \vartheta } 
\prod_{\ell_w \in L_v}(({\Omega^{\prime}}^{-1}_{\nu_{\ell_w}}
{\bf f}_{\alpha_w})\cdot \beta^{\ell_w})\Big\rvert \leq \prod_{m=0}^{\kappa(\alpha)}\left[2\Hat{\Omega}^{-1}(p_{m+1})\right]^{\Hat{N}_k(\vartheta)}A^N B^{|\alpha |}M_{|\alpha |} \, ,
\end{equation}
for some positive constants $A,B$, and $p_{\kappa(\alpha)}\leq |\alpha|< p_{\kappa(\alpha)+1}$. Finally Lemma~\ref{lemma:conteggiovf} gives:
\begin{equation}
\label{eq:finalevf}
\log \frac{|{\bf h}_{\alpha}|}{N_{|\alpha|}} \leq |\alpha|\log C-\log \frac{N_{|\alpha|}}{M_{|\alpha|}}+2|\alpha|\sum_{m=0}^{\kappa(\alpha)}\frac{\log \Omega^{-1}(p_{m+1})}{p_m} \, ,
\end{equation}
for some $C>0$ and the thesis follows dividing by $|\alpha|$ and passing to the limit superior.

\subsection{A result for some analytic vector fields of $\C^2$.}
\label{ssec:anavfc2}

For $2$--dimensional analytic vector field the existence of the continued fraction and the convergents allows us to improve the previous Theorem, giving an optimal (we conjecture) estimate on the ``size'' of the analyticity domain of the linearization.

Mattei and Moussu~\cite{MatteiMoussu} proved, using the holonomy construction, that linearization of germs implies linearization of the foliation~\footnote{This means that a vector field of the form~\eqref{eq:anavf1} can be put in the form:
\begin{equation*}
\begin{cases}
\dot z_1 &= -z_1 (1+h(z_1,z_2) ) \\
\dot z_2 &= \omega z_2(1+h(z_1,z_2)) \, ,         
\end{cases}
\end{equation*}
for some analytic function $h$ such that $h(0,0)=0$.} associated to the vector fields of $(\C^2,0)$. In~\cite{YoccozPerezMarco} authors proved, using Hormander $\bar{\partial}$--techniques, the converse statement. More precisely they proved that the foliation associated to the analytic vector field:
\begin{equation}
\label{eq:anavf1}
\begin{cases}
\dot z_1 &= -z_1 (1+\dots ) \\
\dot z_2 &= \omega z_2(1+\dots) \, ,         
\end{cases}
\end{equation}
where $\omega >0$ and the suspension points mean terms of order bigger than $1$, has the same analytical classification of the germs of $(\C,0)$: $f(z)=e^{2\pi i \omega}z+\bigo{z^2}$. Using the results of~\cite{Yoccoz} they obtain as corollaries that: if $\omega$ is a Bruno number then the foliation associated to~\eqref{eq:anavf1} is analytically linearizable, whereas if $\omega$ is not a Bruno number then there exist analytic vector fields of the form~\eqref{eq:anavf1} whose foliation are not analytically linearizable.

Here we push up this analogy between vector fields and germs, by proving that the linearizing function of the vector field is analytic in domain containing a ball of radius $\rho$ which satisfies the same lower bound (in term of the Bruno function) as the radius of convergence of the linearizing function of the germ~\cite{CarlettiMarmi,Yoccoz} does.

To do this we must introduce some {\em normalization condition} for the vector field; let $\omega >0$, we consider the family $\mathcal{F}_{\omega}$  of analytic vector fields $\mathbf{F}:\mathbb{D}\times\mathbb{D} \rightarrow \mathbb{C}^2$ of the form
\begin{equation}
\begin{cases}
F_1(z_1,z_2) &=-z_1+\sum_{|\alpha|\geq 2}f_{\alpha,1}z^{\alpha} \\
F_1(z_1,z_2) &=\omega z_2+\sum_{|\alpha|\geq 2}f_{\alpha,2}z^{\alpha} \, ,
\end{cases}
\label{eq:vectfieldform}
\end{equation}
with $|f_{\alpha,j}|\leq 1$ for all $|\alpha|\geq 2$ and $j=1,2$.

For power series in several complex variables the analogue of the disk of convergence is the {\em complete Reinhardt domain of center $0$}, $\mathcal{R}_0$, by studying the distance of the origin to the boundary of this domain we can obtain informations about its ``size''. Fixing the non linear part of the vector field: $\mathbf{f}=\sum_{|\alpha|\geq 2}\mathbf{f}_{\alpha}z^{\alpha}$, this distance is given by $d_{\mathbf{F}}=\inf_{(z_1,z_2)\in \mathcal{R}_0}(|z_1|^2+|z_2|^2)^{1/2}$. The family $\mathcal{F}_{\omega}$ is compact w.r.t. the uniform convergence on compact subsets of $\mathbb{D}\times\mathbb{D}$ (use Weierstrass Theorem and Cauchy's estimates in $\C^2$, see for example~\cite{Shabat}) so we can define $d_{\omega}=\inf_{\mathbf{F}\in\mathcal{F}_{\omega}}d_{\mathbf{F}}$.

Let $\rho_{\mathbf{F}}>0$ and let us introduce $P(0,\rho_{\mathbf{F}})=\{ (z_1,z_2)\in \C^2 : |z_i|<\rho_{\mathbf{F}},i=1,2 \}$, the biggest polydisk of center $0$ contained in $\mathcal{R}_0$, whose radius depends on the vector field $\mathbf{F}$. Trivially $\rho_{\mathbf{F}}$ and $d_{\mathbf{F}}$ are related by a coefficient depending only on the dimension: $\sqrt{2}\rho_{\mathbf{F}}=d_{\mathbf{F}}$. We can then prove

\begin{theorem}[Lower bound on $d_{\omega}$]
\label{the:mainvfn2}
Let $\omega >0$ be a Bruno number, then there exists an universal constant $C$ such that:
\begin{equation}
\label{eq:lowboundd}
\log d_{\omega} \geq -B(\omega)-C \, ,
\end{equation}
where $B(\omega)$ is the value of the Bruno function~\cite{MMY} on $\omega$.
\end{theorem}

We don't prove this Theorem being its proof very close to the one of Theorem~\ref{the:mainvf} case 1), we only stress that the use of the continued fraction allows us to give an ``optimal'' counting Lemma as done in~\cite{CarlettiMarmi,Davie}, which essentially bounds the number of lines on scale $k$ in a rooted labeled tree of order $N$ and total momentum $\nu_{\vartheta}$, with $\Big\lfloor \frac{\bar{\nu}_{\vartheta}}{q_k}\Big\rfloor$, being $(q_k)_k$ the denominators of the convergent to $\omega$.

In the case of analytic germs of $(\C,0)$ Yoccoz~\cite{Yoccoz} proved that the same bound holds from above for the radius of convergence of the linearization; the sophisticate techniques used in~\cite{YoccozPerezMarco} would lead to prove: 
\begin{equation*}
\log d_{\omega} \leq -C B(\omega)+C^{\prime} \, ,
\end{equation*}
for some constants $C>1$ and $C^{\prime}$, we {\em conjecture that one can take $C=1$}. We are not able to prove this fact but can prove that the power series obtained replacing the coefficients of the linearization with their absolute values is divergent whenever $\omega$ is not a Bruno number (a similar result has been proved in~\cite{Yoccoz} Appendix $2$ and in~\cite{CarlettiMarmi} paragraph 2.4 for germs).

\begin{remark}[Ultradifferentiable vector fields of $\C^2$]

In the more general case of ultradifferentiable vector fields of $\C^2$ we can improve Theorem~\ref{the:mainvf} showing that we can linearize the vector field under weaker conditions.
\begin{theorem}
Let $\omega >0$ and let $(p_k/q_k)_k$ be its convergents. Let $\mathbf{F}$ be a vector field of the form~\eqref{eq:vectfieldform} (without additional hypotheses on the coefficients $\mathbf{f}_{\alpha}$), let $(M_n)_n$ and $(N_n)_n$ be two sequences verifying conditions 0)--3) of section~\ref{classiMn}. Then
\begin{enumerate}
\item If moreover $\mathbf{F}$ belongs to $\class{M_n}$ and $\omega$ is a Bruno number then also the linearization ${\bf h}$ belongs to $\class{M_k}\cap V_1$.

\item If ${\bf F}$ is a germ of analytic diffeomorphisms of $(\C^2,0)$ and $\omega$ verifies:
\begin{equation}
\limsup_{|\alpha |\rightarrow +\infty} \left(  \sum_{m=0}^{\kappa (\alpha)} \frac{\log q_{m+1}}{q_m}-\frac{1}{|\alpha|}\log N_{|\alpha|} \right) <+\infty \, ,
\end{equation}
where $\kappa(\alpha)$ is the integer defined by: $q_{\kappa(\alpha)}\leq |\alpha| <q_{\kappa(\alpha)+1}$, then the linearization ${\bf h}$ belongs 
to $\class{N_k}\cap V_1$.

\item If ${\bf F}\in \class{M_k}\cap V_0$, the sequence $(M_k)_k$ 
is asymptotically bounded by the sequence $(N_k)_k$ (namely $N_k \geq M_k$ 
for all sufficiently large $k$) and $\omega$ verifies:
\begin{equation}
\limsup_{|\alpha |\rightarrow +\infty} \left( \sum_{m=0}^{\kappa (\alpha)} \frac{\log q_{m+1}}{q_m}-\frac{1}{|\alpha|}\log \frac{N_{|\alpha|}}{M_{|\alpha|}} \right) <+\infty \, ,
\end{equation}
where $\kappa(\alpha)$ is the integer defined by: $q_{\kappa(\alpha)}\leq |\alpha| <q_{\kappa(\alpha)+1}$, then the linearization ${\bf h}$ belongs to $\class{N_k}\cap V_1$.
\end{enumerate}
\end{theorem}

The proof follows closely the one of Theorem~\ref{the:mainvf} and the weaker arithmetical condition are obtained using the ``optimal'' counting function as done in the proof of Theorem~\ref{the:mainvfn2} and in~\cite{CarlettiMarmi,Davie}.
\end{remark}

\appendix

\section{Ultrametric structures and analyticity.}
\label{ultrametricappendix}

Let $k$ be an ultrametric field and let $v$ be a
valuation. The ring $A_v = \{ x \in k \mid v \left( x \right) \geq 0
\} $ is called the {\em ring of the valuation} $v$ and the sets
$I_{\alpha}^{\prime} = \{ x \in A_v \mid v \left( x \right) >
\alpha \} $, for $\alpha \geq 0$, are ideals of $A_v$.
$I_0^{\prime}$ is the maximal ideal of $A_v$ and it is an open set
in the topology induced by the ultrametric absolute value defined
on $k$:
\begin{equation*}
I_0^{\prime}  = \left\{ x \in A_v \mid v \left( x \right) > 0
\right\}  = \left\{ x \in A_v \mid \lvert x \rvert < 1 \right\}
 = B_{0} \left( 0,1 \right),
\end{equation*}
where $B_{0} \left( x,r \right)=\{y\in k\, , \, |x-y|<r\}$. The
field $k^\prime = A_{v}/I_{0}^\prime $ is called the {\em residue
field}.

\indent Let $n$ and $m$ be positive integers, we consider the
set of the formal power series with coefficients in $k^n$ in the
$m$ variables $X_1,\ldots,X_m$, $S_{n,m}=k^n \left[ \left[
X_1,\ldots,X_m \right] \right]$, $F\in S_{n,m}$:
\begin{equation*}
F=\sum_{\alpha\in\mathbb{N}^m}F_{\alpha}X^{\alpha}=
\sum_{\alpha=(\alpha_1,\ldots,\alpha_m)\in\N^m}F_{\alpha}X_1^{\alpha_1}\ldots
X_m^{\alpha_m},
\end{equation*}
with $F_{\alpha}=(F_{\alpha,1},\ldots,F_{\alpha,n})\in
k^n$ for all $\alpha\in\mathbb{N}^m$. We will be
interested in composition problems, so it is natural to set $m=n$
and to define the composition of two elements $F,G$, with
$v(G)\geq 1$, as
\begin{equation}
\label{defcompos} F\circ
G=\sum_{\alpha\in\mathbb{N}^n}F_{\alpha}G^{\alpha}=
\sum_{\alpha_1,\ldots,\alpha_n\in\mathbb{N}}
F_{\alpha_1,\ldots,\alpha_n}(G_{1})^{\alpha_1}\ldots
(G_{n})^{\alpha_n} \, .
\end{equation}
We will set $S_{n,n}=S_n$ and we will introduce some definitions and
properties of $S_n$, but it's clear that they also hold on
$S_{n,m}$ with some small changes.

\indent Let $F =\sum F_{\alpha}X^{\alpha}\in S_n$, then we say
that $F$ {\em converges} in $B(0,r)$, for some $r>0$, if: $\sum_{\alpha\in\mathbb{N}^n}||F_{\alpha}||r^{|\alpha|} <+\infty$, 
where $|\alpha|=\sum_{i=1}^n \alpha_i$. $F$ will be said
convergent in $B_{0}(0,r)$ if it is convergent in $B(0,r^\prime )$
for all $0<r^\prime <r$.

One has the following result (\cite{Serre}, pp. 67--68):

\begin{proposition}
\label{powerseries}
\item{1.} If $F$ is convergent in $B(0,r)$ then there exists $M>0$
such that
\begin{equation}
\label{Cauchy} ||F_{\alpha}||\le Mr^{-|\alpha|} \quad \forall \alpha\in\mathbb{N}^n \, .
\end{equation}
\item{2.} If there exists $M>0$ such that (\ref{Cauchy}) holds for
all $\alpha\in\mathbb{N}^n$, $F$ converges in $B_{0}(0,r)$ and
uniformly in $B(0,r^\prime )$ for all $0<r^\prime <r$.
\item{3.} Let $\tilde{F}\, :\, B_{0}(0,r)\rightarrow k^n$ denote the
continuous function defined as the sum of the series $F\in S_n$
convergent in $B_{0}(0,r)$. Then $\tilde{F}\equiv 0 \iff F = 0$.
\end{proposition}

We can therefore identify a convergent power series $F$ with its
associated function $\tilde{F}$ and vice versa. Let $U$ be an open
set of $k^n$, $\Tilde{G}:U \rightarrow k^n$ is said to be {\em analytic}
in $U$ if for all $x \in U$ there is a formal power series $ G \in
S_{n}$ and a radius $r>0$ such that:
\begin{enumerate}
\item $B_{0} \left( x,r \right) \subset U$,
\item $G$ converges in $B_{0} \left( 0,r \right)$, and for all $y \in B_{0}
\left( 0,r \right)$, $\Tilde{G} \left( x + y \right)= G \left( y
\right)$.
\end{enumerate}
With a slight abuse of notation we will omit the
superscript $\tilde{}$ to distinguish analytic functions from
convergent power series. If $F$ is a convergent power series on $B(0,r)$ we denote
\begin{equation}
\label{norm}
||F||_r=\sup_{\alpha\in\mathbb{N}^n}||F_{\alpha}||r^{|\alpha|}\, ,
\end{equation}
and we define $\mathcal{A}_r(k^n)=\left\{F\in S_n: ||F||_r
<+\infty \right\}$.

\indent Let $U$ be an open set of $k^n$, $x$ a point of $U$, let us consider $G :U \rightarrow V\subset k^m$. A linear function
$L: k^n \rightarrow k^m$ is called a {\em derivative} of $G$ at
$x$ if:

\begin{equation*}
\lim_{\substack{||y || \rightarrow 0 \\ y \neq 0}} \frac {|| G
\left( x + y \right) - G \left( x \right) - Ly ||} {|| y ||}=0.
\end{equation*}
Clearly if the limit exists then the derivative is unique and it
will be denoted by $dG \left( x \right)$. Let $\delta_i=\left( 0,
\dots,1,\dots,0 \right)\in k^n$ the vector with $1$ at the i--th
place, we call
\begin{equation}
D_iG\left( x \right) =dG\left( x \right) \left( \delta_i
\right)\in k^m \label{partialder}
\end{equation}
the i-th partial derivative of $G$ at $x$. Higher order
derivatives are defined analogously.

\indent Let $G = \sum G_{\alpha}X^{\alpha}$ be an element of
$S_n$, then $G_{\alpha} = \frac{D^{\alpha}G(x)}{\alpha !}$, where
for $\alpha=(\alpha_1,\ldots, \alpha_n)\in\mathbb{N}^n$,
$D^{\alpha}=D^{\alpha_1}_1\ldots D^{\alpha_n}_n$ and
$\alpha!=\alpha_1!\ldots \alpha_n!$. Thus $G$ is just the Taylor
series of $G$ at the point $x$.

It is not difficult to prove that any power series $G\in S_{n}$
convergent in $B_{0}(0,r)$ defines an analytic function $G$ in
$B_{0}(0,r)$. However one should be aware of the fact that in
general the local expansion of a function $G$ analytic on $U$ at a
point $x\in U$ such that $B_{0}(x,r)\subset U$ does not
necessarily converge on all of $B_{0}(0,r)$. This is true if one
assumes $k$ to be algebraically closed.

\indent Let $F \in S_{n}$, $F=\sum_{\alpha \in
\N^n}F_{\alpha}X^{\alpha}$,
 let $\beta \in \N^n$ we define the formal derivative of $F$ by
\begin{equation}
\Delta^{\beta} F= \sum_{\alpha \in \N^n,\alpha \geq \beta}
F_{\alpha} \binom{\alpha}{\beta}X^{\alpha-\beta}\, ,
 \label{formalder}
\end{equation}
where we used the notations: $\alpha \geq \beta$ if $\alpha_i \geq \beta_i$ for $1 \leq i \leq n$, $\binom{\alpha}{\beta}=\frac{\alpha!}{\beta!\left( \alpha - \beta \right)!}$, for $\alpha,\beta \in \N^n$ and one can then prove~\cite{Serre}: $\alpha ! \Delta^{\alpha}=D^{\alpha}, \quad  \binom{\alpha +
\beta}{\alpha} \Delta^{\alpha + \beta}= \Delta^{\alpha}
\Delta^{\beta}$.

Finally we note that the composition of two analytic functions is
analytic (\cite{Serre}, p. 70) and that the  following Cauchy
estimates and Taylor formula hold (\cite{HermanYoccoz}, pp.
421-422):

\begin{proposition}
\label{Taylorformulas} Let $r>0$, $s>0$, let
$F\in\mathcal{A}_r(k^n)$, and let $G\, , \, H$ be two elements of
$\mathcal{A}_s(k^n)$, with $||G||_{s}\le r$ and $||H||_{s}\le r$.
Then the following estimates hold:
\begin{itemize}
\item (Cauchy's estimates) $||\Delta^{\alpha} F||_{r} \le
\frac{||F||_{r}}{r^{|\alpha|}}$ for all $\alpha\in\mathbb{N}^n$;
\item $F\circ G \in \mathcal{A}_s(k^n)$ and $||F\circ G||_s\leq
|| F||_r$;
\item (Taylor's formula)
\begin{align*}
 ||F\circ (G+H) - F\circ G ||_s \leq \frac{||F||_{r}}{r} ||H||_{s} \\
 ||F\circ (G+H) - F\circ G - DF\circ G\cdot H||_{s} \le \frac{||F||_{r}}{r^{2}} ||H||_{s}^{2}
\end{align*}
\end{itemize}
\end{proposition}

\section{Proof of the Bruno counting lemma}
\label{sec:brunocl}

In this section we give the proof of Lemma~\ref{lemma:conteggio}: the {\em Bruno counting lemma} for germs. The proof of Lemma~\ref{lemma:conteggiovf} (the vector fields case) can be done in a similar way and we omit it.

\subsection{Proof of Lemma~\ref{lemma:conteggio}}

Let us recall briefly the object of the Lemma. We are considering rooted labeled trees $\vartheta$, any line produces a divisor and we want to count the number of lines producing ``small divisors'', i.e. the number of lines on scale $k$ for some integer $k$. The way these small divisors accumulate give rise to the arithmetical condition needed to prove the convergence of the series involved. We can then prove

\noindent{{\bf Lemma~\ref{lemma:conteggio}}} (Bruno's Counting Lemma). {\em The number of lines on scale $k$ in a rooted labeled tree verifies the following bound:
\begin{equation}
\label{eq:boundcountB}
N_k(\vartheta)\leq 
\begin{cases}
0 \quad &\text{if }|\bar{\nu}_{\vartheta}|< p_k \\
2\Big \lfloor \frac{|\bar{\nu}_{\vartheta}|}{p_k} \Big \rfloor -1\quad &\text{if }|\bar{\nu}_{\vartheta}|\geq p_k \, .
\end{cases}
\end{equation}
where $\lfloor x\rfloor$ denotes the integer part of the real number $x$. We recall that $\nu_{\vartheta}$ is the total momentum of the tree and $\bar{\nu}_{\vartheta}=\nu_{\vartheta}-\beta^{\ell_{v_1}}$, being $\ell_{v_1}$ the root line of $\vartheta$.}

\proof
Consider firstly the case $|\bar{\nu}_{\vartheta}|< p_k$. Because for any $\ell_w\in\vartheta$ we have: $|\nu_{\vartheta}| \geq |\nu_{\ell_w}|$, we conclude that no one line is on scale $k$: $N_k(\vartheta)=0$.

Consider now the case $|\bar{\nu}_{\vartheta}|\geq p_k$. If $\vartheta \in\forest{N}$ for $N\geq 1$ then using the standard decomposition of the tree $\vartheta=(t,\vartheta^1,\dots,\vartheta^t)$, being $t$ the degree of the root $v_1$, $\vartheta^i\in \forest{N_i}$ and $N_1+\dots + N_t = N-1$, we have 
\begin{eqnarray}
N_k(\vartheta) &=s_{\ell_{v_1}}(k)+N_k(\vartheta^1)+\dots + N_k(\vartheta^t) \label{eq:foundrel}\\ \nu_{\vartheta}&=\alpha_{v_1}-\beta_{v_1}+\nu_{\vartheta^1} +\dots +\nu_{\vartheta^t} \, \label{eq:foundrel2}.
\end{eqnarray}
We will prove~\eqref{eq:boundcountB} by induction of the total momentum of the tree. We will distinguish several cases.

\begin{enumerate}
\item[case A) ${\mathbf s_{\ell_{v_1}}(k)=0}$] Because $|\nu_{\vartheta_i}|<|\nu_{\vartheta}|$ for all $i=1,\dots,t$ we can use the induction hypothesis and from~\eqref{eq:foundrel} we get
\begin{align*}
N_k(\vartheta)&= N_k(\vartheta^1)+\dots + N_k(\vartheta^t) \\
               &\leq 2\Big \lfloor \frac{|\bar{\nu}_{\vartheta^1}|}{p_k} \Big \rfloor -1 +\dots +2\Big \lfloor \frac{|\bar{\nu}_{\vartheta^t}|}{p_k} \Big \rfloor -1 \, ,
\end{align*}
from which the thesis follows using~\eqref{eq:foundrel2} and the hypothesis $|\bar{\nu}_{\vartheta}|\geq p_k$.

\item[case B) ${\mathbf s_{\ell_{v_1}}(k)=1}$] We now consider 3 subcases.
\begin{enumerate}
\item[case B.1) $t=0$] Then~\eqref{eq:foundrel} gives $N_k(\vartheta)=1$ and the thesis follows recalling that $|\bar{\nu}_{\vartheta}|\geq p_k$.
\item[case B.2) $t\geq 2$] Then~\eqref{eq:foundrel} gives
\begin{equation*}
N_k(\vartheta)\leq 1+2\Big \lfloor \frac{|\bar{\nu}_{\vartheta^{i_1}}|}{p_k} \Big \rfloor -1+2\Big \lfloor \frac{|\bar{\nu}_{\vartheta^{i_2}}|}{p_k} \Big \rfloor -1 \, ,
\end{equation*}
and again the thesis follows using $|\bar{\nu}_{\vartheta}|\geq p_k$.
\item[case B.3) $t=1$] Equation~\eqref{eq:foundrel} gives $N_k(\vartheta)\leq 2\Big \lfloor \frac{|\bar{\nu}_{\vartheta^{1}}|}{p_k} \Big \rfloor$, so if $\Big \lfloor \frac{|\bar{\nu}_{\vartheta^{1}}|}{p_k} \Big \rfloor < \Big \lfloor \frac{|\bar{\nu}_{\vartheta}|}{p_k} \Big \rfloor$ the thesis follows again. It remains the case $\Big \lfloor \frac{|\bar{\nu}_{\vartheta^{1}}|}{p_k} \Big \rfloor =\Big \lfloor \frac{|\bar{\nu}_{\vartheta}|}{p_k} \Big \rfloor$, namely $|\bar{\nu}_{\vartheta}|-|\bar{\nu}_{\vartheta^{1}}|<p_k$. Let $v_1^{\prime}$ be the root of $\vartheta^1$, let $t^{\prime}$ be its degree and consider the standard decomposition of the subtree $\vartheta^1=(t^{\prime},\vartheta^{1}_1,\dots,\vartheta^{t^{\prime}}_1)$. Lemma~\ref{lemma:brj1} assures that $s_{\ell_{v_1^{\prime}}}(k)=0$, so~\eqref{eq:foundrel} (written for the standard decomposition of $\vartheta^1$) reduces to:
\begin{equation*}
N_k(\vartheta)= 1+N_k(\vartheta^{1}_1)+\dots +N_k(\vartheta^{t^{\prime}}_1) \, .
\end{equation*}
We now consider the cases B.1, B.2 and B.3 for the subtrees $\vartheta^{1}_i$. 
We claim that if case B.1 or B.2 holds the proof is done, whereas in case B.3
 the proof is achieved only if the first subcase holds, but the remaining
 case can happen only a finite number of times, and so in this case too,
 the proof is done.
\end{enumerate}
\end{enumerate}
\endproof

\noindent
{\bf Acknowledgement. }I am are grateful to J.-C. Yoccoz for pointing out the difference between the analytic and non--analytic inversion problems and the relation of the latter with small divisors problems. I also thank A. Albouy for some details about the history of the Lagrange inversion formula, in particular for reference~\cite{Cajori}

\smallskip


\end{document}